\documentclass[11pt,twoside,a4paper]{article}
\usepackage{amsmath,amssymb,mathrsfs}
\usepackage{amsthm}
\usepackage[utf8x]{inputenc}
\usepackage{bm}
\usepackage{avant}
\usepackage[english]{babel}
\usepackage{thmtools,thm-restate}
\usepackage[nottoc]{tocbibind}
\usepackage{hyperref}
\usepackage{graphicx}
\usepackage{enumitem}
\usepackage{tikz}
\usepackage{bbm}
\usepackage[inner=2.4cm,outer=2.4cm,top=2.8cm,bottom=2.8cm]{geometry}
\usepackage{mathtools}
\usepackage{makeidx}
\usepackage{bold-extra}

\DeclareMathOperator{\id}{id}

\DeclareMathOperator{\supp}{supp}
\newcommand{\scal}[2]{\ensuremath{\langle #1 , #2 \rangle}} 
\newcommand{\norm}[1]{\left\lVert#1\right\rVert}
\newcommand{\modu}[1]{\left |#1\right |}
\newcommand{\Leb}{\mathscr{L}}
\newcommand{\setR}{\mathbb{R}}
\newcommand{\setN}{\mathbb{N}}
\newcommand{\R}{\mathbb{R}}
\newcommand{\Q}{\mathbb{Q}}
\newcommand{\p}{\mathtt p} 
\newcommand{\de}{\ensuremath{\, \mathrm d}} 

\newcommand{\suchthat}{\ensuremath{\,:\,}} 
\newcommand\restr[2]{{
  \left.\kern-\nulldelimiterspace 
  #1 
  \right|_{#2} 
  }}

\newcommand{\weakto}{\rightharpoonup}

\newcommand{\CD}{\mathsf{CD}}

\newcommand{\rest}[2]{{{\rm restr}_{#1}^{#2}}}
\newcommand{\di}{\mathsf d} 
\newcommand{\m}{\mathfrak m} 
\DeclareMathOperator{\Ent}{Ent}
\DeclareMathOperator{\Geo}{Geo}
\DeclareMathOperator{\OptGeo}{OptGeo}
\DeclareMathOperator{\OptPlans}{OptPlans}
\newcommand{\Prob}{\mathscr{P}}
\newcommand{\ProbTwo}{\mathscr{P}_2}
\newcommand{\Cost}{\ensuremath{\mathcal C}} 
\newcommand{\ppi}{{\mbox{\boldmath$\pi$}}}

\usepackage{titlesec}
\titleformat{\section}
  {\scshape\large\centering}
  {\thesection}{0.5em}{}
\titleformat{\subsection}
  {\scshape\centering}
  {\thesubsection}{0.4em}{}  

\author{Mattia Magnabosco}
\title{\textbf{A Metric Stability Result for the Very Strict CD Condition}}
\date{}

\newtheoremstyle{remark}
        {10pt}
        {10pt}
        {}
        {}
        {\itshape}
        {.}
        {.4em}
        {}

\newtheoremstyle{proof}
        {10pt}
        {10pt}
        {}
        {}
        {\itshape}
        {.}
        {.4em}
        {}
        
\newtheoremstyle{definition}
        {10pt}
        {10pt}
        {}
        {}
        {\bfseries}
        {.}
        {.4em}
        {}

\newtheoremstyle{theorem}
        {10pt}
        {10pt}
        {\slshape}
        {}
        {\bfseries}
        {.}
        {.4em}
        {}

\theoremstyle{theorem}
\newtheorem{theorem}{Theorem}[section]
\newtheorem{prop}[theorem]{Proposition}
\newtheorem{corollary}[theorem]{Corollary}
\newtheorem{lemma}[theorem]{Lemma}

\theoremstyle{definition}
\newtheorem{definition}[theorem]{Definition}

\theoremstyle{remark}
\newtheorem{remark}[theorem]{Remark}

\theoremstyle{proof}
\newtheorem*{pro}{Proof}
\newenvironment{pr}{\begin{pro}%
 \pushQED{\qed}}%
 {\popQED\end{pro}}

\makeindex
\begin{document}
\maketitle

\begin{abstract}
    In \cite{schultz2017existence} Schultz generalized the work of Rajala and Sturm \cite{rajalasturm}, proving that a weak non-branching condition holds in the more general setting of very strict CD spaces. Anyway, similar to what happens for the strong CD condition, the very strict CD condition seems not to be stable with respect to the measured Gromov Hausdorff convergence (cf. \cite{MM-Example}).\\
    In this article I prove a stability result for the very strict CD condition, assuming some metric requirements on the converging sequence and on the limit space. The proof relies on the notions of \textit{consistent geodesic flow} and \textit{consistent plan selection}, which allow to treat separately the static and the dynamic part of a Wasserstein geodesic. As an application, I prove that the metric measure space $\R^N$ equipped with a crystalline norm and with the Lebesgue measure satisfies the very strict $\CD(0,\infty)$ condition.  
\end{abstract}

In their pivotal works Lott, Villani \cite{lottvillani} and Sturm \cite{sturm2006,sturm2006ii} introduced a weak notion of curvature dimension bounds, which strongly relies on the theory of Optimal Transport. They noticed that, in a Riemannian manifold, a uniform bound on the Ricci tensor is equivalent to the uniform convexity of the Boltzmann-Shannon entropy functional in the Wasserstein space. This allowed them to define a consistent notion of curvature dimension bound for metric measure spaces, that is known as CD condition. The metric measure spaces satisfying the CD condition are called CD spaces and enjoy some remarkable analytic and geometric properties.

The validity of the CD condition in a metric measure space $(X,\di,\m)$ is strongly related to the metric structure of its Wasserstein space, which is in turn strictly dependent on the metric structure of $(X,\di,\m)$. For this reason, it is not surprising that some properties of CD spaces hold under additional metric assumptions. Among them, one of the most important is the non-branching condition, which basically prevents two different geodesic to coincide in an interval of times. Since the first works on CD spaces, it has been clear that the non-branching assumption, associated with the CD condition, could confer some nice properties to a metric measure space. For example, Sturm \cite{sturm2006} was able to prove the tensorization property and the local-to-global property, while Gigli \cite{gigli11} managed to solve the Monge problem. The relation between non-branching assumption and CD condition was made even more interesting by the work of Rajala and Sturm \cite{rajalasturm}. They proved that the strong CD condition implies a weak version of the non-branching one, that they called essentially non-branching. 
The work of Rajala and Sturm was then generalized to the wider context of very strict CD spaces by Schultz in \cite{schultz2017existence} (see also \cite{Schultz2019EquivalentDO} and \cite{Schultz2019OnOO}, where he investigates some properties of very strict CD spaces). In particular, he also underlined that every very strict CD space satisfies a weak non-branching condition, that I will call \textit{weak essentially non-branching}.

One of the most important properties of CD spaces is their stability with respect to the measured Gromov Hausdorff convergence. Unfortunately this rigidity result cannot hold for the strong CD condition and, accordingly to \cite{MM-Example}, it also does not hold for the so called \textit{strict} CD condition, which is (a priori) weaker than the very strict one, but stronger than the weak one.
In particular, as explained in \cite{MM-Example}, it is not possible to deduce in general any type non-branching condition for a measured Gromov Hausdorff limit space. This motivates to add some analytic or metric assumption on the converging spaces, in order to achieve non-branching at the limit. In this direction the most remarkable result is provided by the theory of RCD spaces, pioneered by Ambrosio Gigli and Savaré in \cite{Ambrosio_2013} and \cite{Ambrosio_2014}. In fact these spaces are stable with respect to the measured Gromov Hausdorff convergence and essentially non-branching. In this article I present a stability result for very strict CD spaces, assuming metric requirements on the converging sequence and on the limit space. 

In particular, the first section is dedicated to introduce the necessary preliminary notions, related both to the Optimal Transport theory and to CD conditions. In this sense, this section should be understood as a list of prerequisites and not as a complete treatment of the basic theory. For a full and precise discussions about the Optimal Transport theory I refer the reader \cite{AmbrosioNotes2}, \cite{ambrosio2013user}, \cite{villani2003} and \cite{villani2008}.

In the second section I prove a purely metric stability result, which assume some strong rigidity requirements, but nevertheless can be applied to a fair variety of metric measure spaces. This result relies on the notions of \textit{consistent geodesic flow} and \textit{consistent plan selection}, which, as it will be clear in the following, allow me to treat separately the dynamic and the static parts of Wasserstein geodesics. The proof of this result uses an approximation strategy, and it is completely different from the arguments used for the RCD spaces theory.

The result of the second section can be applied to the metric measure space $\R^N$ equipped with a crystalline norm and with the Lebesgue measure, this is explained in the last section. In particular I will show how a secondary variational problem can provide a consistent plan selection, associated to the Euclidean consistent geodesic flow. This will allow to conclude that these metric measure spaces are very strict CD spaces, and therefore they are weakly essentially non-branching.

\section{Preliminary Notions}
This first section is aimed to state all the preliminary results I will need in the following.

\subsection{The Optimal Transport Problem}

The original formulation of the Optimal Transport problem, due to Monge, dates back to the XVIII century, and it is the following: given two topological spaces $X,Y$, two probability measures $\mu\in\Prob(X)$, $\nu\in\Prob(Y)$ and a non-negative Borel cost function $c:X\times Y \to  [0,\infty]$, look for the maps $T$ that minimize the following quantity 
\begin{equation} \label{MongeFormulation}
	\inf\left\{ \int_X c(x,T(x)) \de \mu(x) \suchthat \text{$T:X\to Y$ Borel, $T_\# \mu =\nu$} \right\}. \tag{M}
\end{equation}
The minimizers of the Monge problem are called optimal transport maps and in general do not necessarily exist. Therefore for the development of a general theory, it is necessary to introduce a slight generalization, due to Kantorovich. Defining the set of transport plans from $\mu$ to $\nu$:
\begin{equation*}
\Gamma(\mu,\nu) := \{ \pi\in\Prob(X\times Y) \suchthat (\p_X)_\#\pi = \mu \,\,\text{and} \,\, (\p_Y)_\#\pi = \nu\},
\end{equation*}
the Kantorovich's formulation of the optimal transport problem asks to find minima and minimizers of
\begin{equation} \label{KantorovichFormulation}
\Cost(\mu,\nu):=	\inf \left\{ \int_{X\times Y} c(x,y) \de \pi(x,y) \suchthat \pi\in\Gamma(\mu,\nu) \right\}. \tag{K}
\end{equation}
This new problem admits minimizers under weak assumptions, in fact the following theorem holds.
\begin{theorem}[Kantorovich] \label{thm:MinimumInKantorovich}
	Let $X$ and $Y$ be Polish spaces and $c:X\times Y \to [0,\infty]$ a lower semicontinuous cost function, then the minimum in the Kantorovich's formulation \eqref{KantorovichFormulation} is attained.
\end{theorem}
The minimizers of the Kantorovich problem are called optimal transport plans and the set of all of them will be denoted by $\OptPlans(\mu,\nu)$. Notice that this set obviously depends on the cost function $c$, anyway I will usually avoid to make this dependence explicit, since many times it will be clear from the context. A transport plan $\pi\in\Gamma(\mu,\nu) $ is said to be induced by a map if there exists a $\mu$-measurable map $T:X \to Y$ such that $\pi=(\id,T)_\# \mu$. Notice that these transport plans are exactly the ones considered in the Monge's minimization problem \eqref{MongeFormulation}.

\begin{remark}\label{rem:uniquenessthroughmap}
Suppose that every minimizer of the Kantorovich problem between the measures $\mu,\nu\in\Prob(X)$ is induced by a map, and thus is a minimizer for the
Monge problem. Then the Kantorovich problem between $\mu$ and $\nu$ admit a unique minimizer, which is clearly induced by a map. In fact, given two distinct transport plans $\pi_1=(\id, T_1)_\#\mu,\,\pi_2=(\id, T_2)_\#\mu\in \OptPlans(\mu,\nu)$, their combination $\pi=\frac{\pi_1+\pi_2}{2}$ is itself an optimal plan and it is not induced by a map, contradicting the assumption.
\end{remark}

\noindent A fundamental approach in facing the Optimal Transport problem is the one of $c$-duality, which allows to prove some very interesting and useful results. Before stating them let me introduce the notions of $c$-cyclical monotonicity, $c$-conjugate function and $c$-concave function.

\begin{definition}
A set $\Gamma\subset X\times Y$ is said to be $c$-cyclically monotone if 
\begin{equation*}
    \sum_{i=1}^{N} c\left(x_{i}, y_{\sigma(i)}\right) \geq \sum_{i=1}^{N} c\left(x_{i}, y_{i}\right)
\end{equation*}
for every $N\geq1$, every permutation $\sigma$ of $\{1,\dots,N\}$ and every $(x_i,y_i)\in \Gamma$ for $i=1,\dots,N$.
\end{definition}

\begin{definition}
Given a function $\phi:X \to \R \cup\{-\infty\}$, define its $c$-conjugate function $\phi^c$ as 
\begin{equation*}
    \phi^{c}(y):=\inf _{x \in X}\{c(x, y)-\phi(x)\}.
\end{equation*}
Analogously, given $\psi:Y \to \R \cup\{-\infty\}$, define its $c$-conjugate function $\psi^c$ as 
\begin{equation*}
    \psi^{c}(x):=\inf _{y \in Y}\{c(x, y)-\psi(y)\}.
\end{equation*}
\end{definition}

\noindent Notice that, by definition, given a function $\phi:X \to \R \cup\{-\infty\}$, $\phi(x)+\phi^c(y)\leq c(x,y)$ for every $(x,y)\in X \times Y$. 

\begin{definition}
A function $\phi:X \to \R \cup\{-\infty\}$ is said to be $c$-concave if it is the infimum of a family of $c$-affine functions $c(\cdot, y) + \alpha$. Analogously, $\psi:Y \to \R \cup\{-\infty\}$ is said to be $c$-concave if it is the infimum of a family of $c$-affine functions $c(x, \cdot) + \beta$.
\end{definition}

\noindent The first important result of the $c$-duality theory is the following, which summarize the strict relation that there is between optimality and $c$-cyclical monotonicity.

\begin{prop}\label{prop:cmonotonicity}
Let $X$ and $Y$ be Polish spaces and $c:X\times Y \to [0,\infty]$ a lower semicontinuous cost function. Then every optimal transport plan $\pi\in \OptPlans(\mu,\nu)$ such that $\int c \de \pi<\infty$ is concentrated in a $c$-cyclically monotone set. 
Moreover, if there exist two functions $a\in L^1(X,\mu)$ and $b\in L^1(Y,\nu)$ such that $c(x,y)\leq a(x)+b(x)$ for every $(x,y)\in X\times Y$, a plan $\pi\in \Gamma(\mu,\nu)$ is optimal only if it is concentrated on a $c$-cyclically monotone set.
\end{prop}

\noindent The $c$-duality theory also allows to state the following duality result. 

\begin{prop}\label{prop:Kantorovichpotential}
 Let $X$ and $Y$ be Polish spaces and $c:X\times Y \to [0,\infty]$ a lower semicontinuous cost function. If there exist two functions $a\in L^1(X,\mu)$ and $b\in L^1(Y,\nu)$ such that $c(x,y)\leq a(x)+b(x)$ for every $(x,y)\in X\times Y$, then there exists a $c$-concave function $\phi:X \to \R \cup \{-\infty\}$ satisfying 
 \begin{equation*}
     \Cost(\mu,\nu)= \int \phi \de \mu + \int \phi^c \de \nu. 
 \end{equation*}
 Such a function $\phi$ is called Kantorovich potential.
\end{prop}

\begin{remark}\label{remark:Kantorovich}
Notice that, if the cost $c$ is continuous, every $c$-concave function is upper semicontinuous, being the infimum of continuous functions. Therefore, according to Proposition \ref{prop:Kantorovichpotential}, it is possible to find an upper semicontinuous Kantorovich potential $\phi$ and its $c$-conjugate function $\phi^c$ will be itself upper semicontinuous.
\end{remark}

\subsection{The Wasserstein Space and the Entropy Functional}

In this section I am going to consider the Optimal Transport problem in the special case in which $X=Y$, $(X,\di)$ is a Polish metric space and the cost function is equal to the distance squared, that is $c(x,y)=\di^2(x,y)$. In this context the Kantorovich's minimization problem induces the so called Wasserstein distance on the space $\ProbTwo(X)$ of probabilities with finite second order moment. Let me now give the precise definitions.

\begin{definition}
Define the set
\begin{equation*}
    \ProbTwo(X):= \left\{\mu\in \Prob(X) \suchthat \int \di^2(x,x_0) \de\mu(x) <\infty \text{ for one (and thus all) } x_0\in X\right\}
\end{equation*}
\end{definition}

\begin{definition}[Wasserstein distance]
Given two measures $\mu,\nu\in \ProbTwo(X)$ define their Wasserstein distance $W_2(\mu,\nu)$ as 
\begin{equation*}
W_2^2(\mu,\nu) := \min \left\{ \int d^2(x,y) \de \pi(x,y) \suchthat \pi \in \Gamma(\mu,\nu) \right \}.
\end{equation*}
\end{definition}

\begin{prop}
 $W_2$ is actually a distance on $\ProbTwo(X)$ and $(\ProbTwo(X),W_2)$ is a Polish metric space.
\end{prop}

\noindent The convergence of measures in $\ProbTwo(X)$ with respect to the Wasserstein distance can be easily characterized and this is very useful in many situation. 

\begin{prop}\label{prop:chrW2covergence}
Let $(\mu_n)_{n\in \setN} \subset \ProbTwo(X)$ be a sequence of measures and let $\mu\in \ProbTwo(X)$, then $ \mu_n \xrightarrow{W_2} \mu $ if and only if $\mu_n \weakto \mu $ and
 \begin{equation*}
  \int \di^2(x,x_0) \de \mu_n \to \int \di^2(x,x_0) \de \mu \qquad \text{for every }x_0\in X.
 \end{equation*}
 In particular, if $(X,\di)$ is a compact metric space, the Wasserstein convergence is equivalent to weak convergence.
\end{prop}

Let me now deal with the geodesic structure of $(\ProbTwo(X),W_2)$, which, as the following statement shows, is heavily related to the one of the base space $(X,\di)$. This fact makes the Wasserstein space very important, and allows to prove many remarkable facts. First of all, notice that every measure $\ppi\in\Prob(C([0,1],X))$ induces a curve $[0,1]\ni t \to \mu_t=(e_t)_\# \ppi \in \Prob(X)$, therefore in the following I will consider measures in $\Prob(C([0,1],X))$ in order to consider curves in the Wasserstein space.

\begin{prop}\label{prop:optgeo}
If $(X,\di)$ is a geodesic space then $(\ProbTwo(X),W_2)$ is geodesic as well.
More precisely, given two measures $\mu,\nu\in \ProbTwo(X)$, the measure $\ppi\in\Prob(C([0,1],X))$ is a constant speed Wassertein geodesic connecting $\mu$ and $\nu$ if and only if it is concentrated in $\Geo(X)$ (that is the space of constant speed geodesics in $(X,\di)$) and $(e_0,e_1)_\#\ppi\in\OptPlans(\mu,\nu)$. In this case it is said that $\ppi$ is an optimal geodesic plan between $\mu$ and $\nu$ and this will be denoted as $\ppi\in\OptGeo(\mu,\nu)$.
\end{prop}

\begin{remark}\label{rem:uniquenessthroughmap2}
I will say that an optimal geodesic plan $\ppi\in\OptGeo(\mu,\nu)$ is induced by a map if there exists a $\mu$-measurable map $\Theta:X\to \Geo(X)$, such that $\ppi= \Theta_\# \mu$.
Following the argument explained in Remark \ref{rem:uniquenessthroughmap}, it is possible to conclude that, if any optimal geodesic plan $\ppi\in\OptGeo(\mu,\nu)$ between two given measures $\mu,\nu\in\ProbTwo(X)$ is induced by a map, then there exists a unique $\ppi\in\OptGeo(\mu,\nu)$ and it is obviously induced by a map.
\end{remark}

Let me now introduce the entropy functional that will be the main character in defining the notion of weak curvature dimension bounds. As it will be soon clear, the most appropriate framework in which deal with the entropy functional, is the one of metric measure spaces.

\begin{definition}
A metric measure space is a triple $(X,\mathsf{d},\mathfrak{m})$, where $(X,\mathsf{d})$ is a Polish metric space and $\mathfrak{m}$ is a non-negative and non-null Borel measure on $X$, finite on bounded sets. Moreover, a quadruple $(X,\di,\m,p)$ is called pointed metric measure space if $(X,\di,\m)$ is a metric measure space and $p$ is a point in $X$.
\end{definition}

\begin{remark}
In this article I assume that every metric measure space I am going to consider satisfies the following estimate
\begin{equation}\label{eq:growthcondition}
    \int e^{-c \cdot \di^2(x,x_0)} \de \m(x) <\infty,
\end{equation}
for some (and thus all) $x_0\in X$ and a suitable constant $c>0$. This is essentially a technical assumption, but it is useful to ensure the lower semicontinuity of the entropy functional (see Proposition \ref{prop:lscofent}).
\end{remark}

\noindent Let me now properly define the entropy functional.

\begin{definition}
In a metric measure space $(X,\di,\m)$, given a measure $\nu\in \mathcal{M}(X)$ define the relative entropy functional with respect to $\nu$ $\Ent_\nu:\ProbTwo(X)\to\overline\R $ as 
\begin{equation*}
    \Ent_\nu(\mu):= 
    \begin{cases}
    \int \rho \log \rho \de \nu  &\text{if }\mu\ll\nu \text{ and } \mu=\rho\nu\\
    +\infty &\text{otherwise}
    \end{cases}.
\end{equation*}
The entropy functional relative to the reference measure $\m$ will be simply denoted by $\Ent$. 
\end{definition}

\begin{remark}
According to \cite{Ambrosio_2015}, condition \eqref{eq:growthcondition} prevents the entropy functional $\Ent$ to take the value $-\infty$.
\end{remark}

\noindent The most important property of the entropy functional is its lower semicontinuity with respect to the different notions of convergence in spaces of probabilities. Some results in this direction are collected in this proposition.

\begin{prop}\label{prop:lscofent}
 If $\m(X)<\infty$ the functional $\Ent$ is lower semicontinuous with respect to the weak topology of $\ProbTwo(X)$, while if $\m(X)=\infty$ (but \eqref{eq:growthcondition} holds) $\Ent$ is lower semicontinuous with respect to the Wasserstein convergence.
\end{prop}

I conclude this subsection introducing two more spaces of probabilities, that will play an important role in the following. 

\begin{definition}
Introduce the space $\Prob_{ac}(X)\subseteq\ProbTwo(X)$ of probabilities absolutely continuous with respect to $\m$, with finite second order moments. Define also the space $\Prob^*(X)\subseteq\Prob_{ac}(X)$ as 
\begin{equation*}
    \Prob^*(X) := \{\mu\in\ProbTwo(X)\suchthat\Ent(\mu)<\infty\}.
\end{equation*}
\end{definition}

\subsection{Curvature Dimension Bounds}

In this section I introduce the notions of curvature dimension bound and CD space, presenting also the results which are the starting point of this work. Let me begin with the definition of weak and strong CD condition.

\begin{definition}\label{def:cdinfty}
A metric measure space $(X,\di,\m)$ is said to satisfy the (weak) $\CD(K,\infty)$ condition and to be a (weak) $\CD(K,\infty)$ space, if for every $\mu_0,\mu_1\in\Prob^*(X)$ there exists a Wasserstein geodesic with constant speed $(\mu_t)_{t\in[0,1]}\subset\Prob^*(X)$ connecting them, along which the relative entropy functional is $K$-convex, that is 
\begin{equation*}
    \Ent(\mu_t) \leq (1-t) \Ent(\mu_0) + t\Ent(\mu_1) -t(1-t)\frac K2  W_2^2(\mu_0,\mu_1), \qquad \text{for every }t\in [0,1].
\end{equation*}
Moreover $(X,\mathsf{d},\mathfrak{m})$ is said to satisfy a strong $\CD(K,\infty)$ condition and to be a strong $\CD(K,\infty)$ space if, for every $\mu_0,\mu_1\in\Prob^*(X)$, the relative entropy functional is $K$-convex along every Wasserstein geodesic with constant speed connecting them. 
\end{definition}

\noindent The following proposition due to Villani \cite{villani2008} ensures the validity of CD condition in some familiar metric measure spaces and it will be fundamental in the last section.

\begin{prop}\label{prop:CDinR^n}
Let $\norm{\cdot}$ be a norm on $\R^n$ and let $\di$ be the associated distance, then the metric measure space $(\R^n, \di, \Leb^n)$ is a (weak) $\CD(0,\infty)$ space.
\end{prop}

The next result states the stability of CD condition with respect to the (pointed) measured Gromov Hausdorff convergence. I am not interested in giving a precise the definition of this notion of convergence, because in this article I will only deal with a different and stronger convergence for metric measure spaces. For a precise definition I refer the reader to \cite{villani2008}, where also the next theorem is proven. Let me also point out that the measured Gromov Hausdorff convergence can be in some sense metrized by the $\mathbb{D}$ distance, introduced by Sturm in \cite{sturm2006ii}. Moreover in \cite{Gigli_2015} Gigli, Mondino and Savaré showed that some different notion of convergence for pointed metric measure spaces are equivalent to the pointed measured Gromov Hausdorff convergence.

\begin{theorem}\label{thm:stabilitynoncompact}
Let $(X_k,\di_k,\m_k,p_k)_{k\in\setN}$ be a sequence of locally compact pointed metric measure spaces converging in the pointed measured Gromov Hausdorff sense to a locally compact pointed metric measure space $(X,\di,\m,p)$. Given $K\in\R$, if each $(X_k,\di_k,\m_k)_{k\in\setN}$ satisfies the weak curvature dimension condition $\CD(K,\infty)$, then also $(X,\di,\m)$ satisfies $\CD(K,\infty)$.
\end{theorem}

I am now going to present the Rajala-Sturm theorem, which is the starting point of this work. In order to do this I have to preliminary introduce the notion of essentially non-branching metric measure space.

\begin{definition}
 A metric measure space $(X,\mathsf{d},\mathfrak{m})$ is said to be essentially non-branching if for every absolutely continuous measures $\mu_0,\mu_1\in\Prob_{ac}(X)$, every optimal geodesic plan $\eta$ connecting them is concentrated on a non-branching set of geodesics.
\end{definition}

\begin{theorem}\label{thm:RajalaSturm}
Every strong $\CD(K,\infty)$ metric measure space is essentially non-branching.
\end{theorem} 

\noindent The work of Rajala and Sturm was then generalized by Schultz \cite{schultz2017existence} and applied to the context of very strict CD spaces.

\begin{definition}\label{def:verystrict}
A metric measure space $(X,\mathsf{d},\mathfrak{m})$ is called a very strict $\CD(K,\infty)$ space
if for every absolutely continuous measures $\mu_0,\mu_1\in\Prob_{ac}(X)$
there exists an optimal geodesic plan $\eta\in\OptGeo(\mu_0,\mu_1)$, so that the entropy functional $\Ent$ satisfies the K-convexity inequality along $(\operatorname{restr}_{t_0}^{t_1})_\# (f\eta)$
for every $t_0<t_1\in [0,1]$, and for all bounded Borel functions $f : \Geo(X) \to \setR^+$ with $\int f \de \eta=1$.
\end{definition}

This CD condition is intermediate between the weak and the strong one and it easy to realize that it cannot imply the essentially non-branching property. Anyway, as pointed out by Schultz, it is possible to prove a weaker version of the non-branching condition. 

\begin{definition}[Weak Essentially Non-Branching]
A metric measure space $(X,\mathsf{d},\mathfrak{m})$ is said to be weakly essentially non-branching if for every absolutely continuous measures $\mu_0,\mu_1\in\ProbTwo(X)$, there exists an optimal geodesic plan connecting them, that is concentrated on a non-branching set of geodesics.
\end{definition}

\begin{theorem}\label{thm:schultz}
Every very strict $\CD(K,\infty)$ space is weakly essentially non-branching.
\end{theorem}

Unfortunately, as the reader can easily notice, the strong CD condition is not stable with respect to the measured Gromov Hausdorff convergence. Moreover, the results in \cite{MM-Example} suggest that it is not possible to prove a general stability result also for the very strict CD condition. These observations motivate this article, where I assume some metric requirements on the converging sequence and on the limit space, in order to prove the very strict CD condition for suitable measured Gromov Hausdorff limit spaces.

\section{A Metric Stability Result}\label{section:ametricstabilityresult}

In this section I state and prove some results that allow to prove very strict condition, and thus weak essentially non-branching, for some special measured Gromov Hausdorff limit spaces. These results do not assume any analytic requirement and are purely metric, therefore they can be applied to a wide variety of metric measure spaces. The way to prove non-branching at the limit in this case is very different from the one used by Ambrosio, Gigli and Savaré in \cite{Ambrosio_2014} and it is actually more straightforward.

First of all, let me introduce two notions which provide a nice strategy to prove the very strict CD condition, they are called consistent geodesic flow and consistent plan selection. As it will be clear in the proof of Theorem \ref{thm:verystrict}, these two concepts allow to decouple the static part from the dynamic one, taking full advantage of Proposition \ref{prop:optgeo}. This, associated with suitable assumption, makes easier to deal with restriction of optimal geodesic plans and thus to prove the very strict CD condition.

\begin{definition}\label{def:nicegeodesicflow}
Let $(X,\di)$ be a metric space. A measurable map $G:X\times X\to C([0,1],X)$ is called consistent geodesic flow if the following properties hold:
\begin{itemize}
    \item[1)] for every $x,y\in X$, $G(x,y)$ is a constant speed geodesic connecting $x$ and $y$, that is $G(x,y)\in \Geo (X)$ with $G(x,y)(0)=x$ and $G(x,y)(1)=y$,  
    \item[2)] $\rest{s}{t} G(x,y) = G\big(G(x,y)(s),G(x,y)(t)\big)$ for every $s<t\in (0,1)$ and every $x,y\in X$.
\end{itemize}
A consistent geodesic flow $G$ is said to be $L$-Lipschitz if 
\begin{equation*}
  \sup_{t\in [0,1]}   \di \big(G(x_1,y_1)(t), G(x_2,y_2)(t)\big)  \leq L \cdot \big( \di^2 (x_1,x_2) + \di^2 (y_1,y_2)\big)^\frac12,
\end{equation*}
 i.e. if it is an $L$-Lipschitz map considered as 
 \begin{equation*}
     G: (X\times X, \di_2) \to \big(\Geo (X), \norm{\cdot}_\text{sup} \big),
 \end{equation*}
where $\di_2=\di \otimes \di$.
\end{definition}

\begin{definition}\label{def:planselection}
Let $(X,\di,\m)$ be a metric measure space and assume there exists a consistent geodesic selection $G$ for the metric space $(X,\di)$. A map $\Pi:\Prob_{ac}(X) \times \Prob_{ac}(X) \to \Prob(X\times X)$ is called consistent plan selection, associated to the flow $G$ if 
\begin{itemize}
    \item[1)] $\Pi(\mu,\nu)\in \OptPlans(\mu,\nu)$ for every $\mu,\nu\in\Prob_{ac}(X)$
    \item[2)] For every $\mu,\nu\in \Prob_{ac}(X)$, every pair of times $s<t\in [0,1]$ and every bounded Borel function $f:X\times X\to \R^+$ with $\int f \de \Pi(\mu,\nu)=1$, if 
    \begin{equation*}
        (G_s)_\# \big(f \cdot \Pi(\mu,\nu)\big), (G_t)_\# \big(f \cdot \Pi(\mu,\nu)\big) \in \Prob_{ac}(X),
    \end{equation*}
    where $G_r$ denotes the map $e_r \circ G$ for every $r\in[0,1]$, then it holds
    \begin{equation*}
        (G_s,G_t)_\# \big(f \cdot \Pi(\mu,\nu)\big)= \Pi\big((G_s)_\# \big(f\cdot \Pi(\mu,\nu)\big), (G_t)_\# \big(f\cdot\Pi(\mu,\nu)\big)\big).
    \end{equation*}
\end{itemize}
\end{definition}

\noindent Before going on I present the following lemma, that provides a useful equivalent characterization for condition 2 in the last definition.

\begin{lemma}\label{lem:conditions}
Condition 2 in Definition \ref{def:planselection} is equivalent to the combination of the following two requirements 
\begin{itemize}
    \item[2.1)] $f \cdot \Pi(\mu,\nu) = \Pi\big( (\p_1)_\# \big(f\cdot \Pi(\mu,\nu)\big),(\p_2)_\# \big(f\cdot \Pi(\mu,\nu)\big)\big)$  for every $\mu,\nu\in \Prob_{ac}(X)$ and every bounded Borel function $f:X\times X\to \R^+$ with $\int f \de \Pi(\mu,\nu)=1$.
    \item[2.2)] For every $\mu,\nu\in \Prob_{ac}(X)$ and every $s<t\in [0,1]$, if
    \begin{equation*}
        (G_s)_\# \Pi(\mu,\nu), (G_t)_\# \Pi(\mu,\nu) \in \Prob_{ac}(X),
    \end{equation*}
    then it holds
    \begin{equation*}
        (G_s,G_t)_\# \Pi(\mu,\nu)= \Pi\big((G_s)_\# \Pi(\mu,\nu),(G_t)_\# \Pi(\mu,\nu)\big).
    \end{equation*}
\end{itemize}
\end{lemma}

\begin{pr}
First of all notice that, putting $f\equiv 1$ in condition 2, one obtains condition 2.2. Moreover, also condition 2.1 can be deduced by condition 2, considering only the case where $s=0$ and $t=1$. Therefore condition 2 implies the combination of 2.1 and 2.2.

On the other hand, assume that both 2.1 and 2.2 hold, then for every $s<t\in [0,1]$, every $\mu,\nu\in \Prob_{ac}(X)$ and every bounded Borel function $f:X\times X\to \R^+$ with $\int f \de \Pi(\mu,\nu)=1$, if 
    \begin{equation*}
        (G_s)_\# \big(f \cdot \Pi(\mu,\nu)\big), (G_t)_\# \big(f \cdot \Pi(\mu,\nu)\big) \in \Prob_{ac}(X),
    \end{equation*}
    it holds that
 \begin{align*}
     (G_s,G_t)_\# (f\cdot \pi) &=  (G_s,G_t)_\# \Pi\big( (\p_1)_\#(f\cdot \pi), (\p_2)_\#(f\cdot \pi)\big)\\
     &= \Pi\bigg( (G_s)_\# \Pi\big( (\p_1)_\#(f\cdot \pi), (\p_2)_\#(f\cdot \pi)\big), (G_t)_\# \Pi\big( (\p_1)_\#(f\cdot \pi), (\p_2)_\#(f\cdot \pi)\big)\bigg) \\
     &= \Pi \bigg((G_s)_\# (f\cdot \pi), (G_t)_\# (f \cdot \pi) \bigg),
 \end{align*}
where I have denoted by $\pi$ the plan $\Pi(\mu,\nu)$, in order to ease the notation. This last relation is exactly the requirement of condition 2. 
\end{pr}

I have introduced everything I need to prove one of the crucial results of this section. It shows how the existence of a consistent geodesic flow and a consistent plan selection, satisfying suitable assumptions, ensures the validity of the very strict CD condition.

\begin{theorem}\label{thm:verystrict}
Given a metric measure space $(X,\di, \m)$, assume there exist a consistent geodesic flow $G$ for $(X,\di)$ and a consistent plan selection $\Pi$ associated to $G$. Suppose also that for every pair of measures $\mu,\nu\in \Prob_{ac}(X)$, the $K$-convexity inequality of the entropy is satisfied along the Wasserstein geodesic $G_\# \Pi(\mu,\nu)$ for a suitable $K$, that is 
\begin{equation*}
    \Ent\big((G_t)_\# \Pi(\mu,\nu) \big) \leq (1-t) \Ent(\mu)+ t\Ent(\nu) - \frac K 2 t(1-t) W_2^2(\mu,\nu),
\end{equation*}
for every $t\in(0,1)$. Then $(X,\di,\m)$ is a very strict $\CD(K,\infty)$ space.
\end{theorem}

\begin{pr}
Fix two measures $\mu,\nu\in \Prob_{ac}(X)$ and call $\pi= \Pi(\mu,\nu)$. Then I need to prove that  the $K$-convexity inequality of the entropy holds along the optimal geodesic plan $(\rest{s}{t})_\# \big( f\cdot G_\#  \pi\big)$, for every $s<t \in [0,1]$ and every bounded Borel function $f:C([0,1],X)\to \R^+$ with $\int f \de \pi=1$. This is obviously true when at least one of its marginals at time $0$ and $1$ is not absolutely continuous, therefore I can assume that
\begin{equation}\label{eq:absolutecont}
    (e_s)_\# \big( f \cdot G_\# \pi\big), (e_t)_\# \big( f \cdot G_\# \pi\big) \in \Prob_{ac}(X).
\end{equation}
In particular this allows me to apply condition 2 in Definition \ref{def:planselection}. Now notice that, since $G$ is obviously injective, if one calls $\tilde f= f \circ G$ it holds
\begin{equation*}
    (\rest{s}{t})_\# \big( f\cdot G_\#  \pi\big) = (\rest{s}{t})_\# \big( G_\# (\tilde f \cdot \pi)\big) = (\rest{s}{t} \circ  G)_\# (\tilde f \cdot \pi).
\end{equation*}
Observe now that the definition of consistent geodesic flow ensures that $\rest{s}{t} \circ  G = G \circ (G_s,G_t)$, thus  
\begin{align*}
    (\rest{s}{t})_\# \big( f\cdot G_\#  \pi\big) = (G \circ (G_s,G_t))_\# (\tilde f\cdot \pi)&= G_\# \big( (G_s,G_t)_\# (\tilde f \cdot \pi)\big) \\
    &= G_\# \Pi\big( (G_s)_\# (\tilde f \cdot \pi), (G_t)_\# (\tilde f \cdot \pi)\big) .
\end{align*}
On the other hand it is obvious that 
\begin{equation*}
    (G_s)_\# (\tilde f \cdot \pi) = (e_s)_\# \big( G_\# (\tilde f \cdot \pi) \big) = (e_s)_\# \big( f \cdot G_\# \pi\big),
\end{equation*}
and similarly $(G_t)_\# (\tilde f \cdot \pi) = (e_t)_\# \big( f \cdot G_\# \pi\big)$, so I can conclude that
\begin{equation*}
    (\rest{s}{t})_\# \big( f\cdot G_\#  \pi\big) = G_\# \Pi\big((e_s)_\# \big( f \cdot G_\# \pi\big), (e_t)_\# \big( f \cdot G_\# \pi\big) \big).
\end{equation*}
At this point the fact that the entropy functional is $K$-convex along $(\rest{s}{t})_\# \big( f\cdot G_\#  \pi\big)$ is an easy consequence of the assumption of the theorem, associated with \eqref{eq:absolutecont}.
\end{pr}

In the remaining part of the section I show how Theorem \ref{thm:verystrict} can be applied in order to prove the very strict CD condition for some suitable measured Gromov Hausdorff limit spaces. The first result I want to present provides sufficient conditions to ensure the existence of a consistent geodesic flow for a limit space. The reader must notice that I am considering a notion of convergence, that is much stronger than the measured Gromov Hausdorff convergence. This choice allows me to avoid some tedious technical details, but it is easy to notice that this result can be somehow extended to measured Gromov Hausdorff limit spaces. Anyway, as the next section confirms, in many easy applications this stronger hypothesis is sufficient.  

\begin{prop}\label{prop:ngf}
Let $(X,\di, \m)$ be a compact metric measure space and let $\{\di_n\}_{n\in \setN}$ be a sequence of distances on $X$ (inducing the same topology) such that there exist a sequence $\{\varepsilon_n\}_{n\in \setN}$ converging to zero satisfying  
\begin{equation*}
   | \di_n(x,y)-\di(x,y)| < \varepsilon_n \quad \text{for every } x,y\in X,
\end{equation*}
in particular the sequence $(X,\di_n,\m)$ measured Gromov Hausdorff converges to $(X,\di, \m)$ by means of the identity map. Given $L>0$, assume that for every $n$ there exists an $L$-Lipschitz consistent geodesic flow $G_n$ for the metric measure space $(X,\di_n, \m)$. Then there exists an $L$-Lipschitz consistent geodesic flow $G$ for the metric measure space $(X,\di, \m)$ and, up to subsequences, $G_n$ converges uniformly to $G$. 
\end{prop}

\begin{pr}
The space $X\times X$ is compact and thus separable, therefore there exists a countable dense set $D\subseteq X\times X$. Notice that for every $(x,y)\in X \times X$ and every $t\in[0,1]$, the sequence $\{G_n(x,y)(t)\}_{n\in\setN}$ is contained in the compact set $X$. Then the diagonal argument ensures that, up to taking a suitable subsequence, there exists a function
\begin{equation*}
    G:D \to \{f:[0,1]\cap \Q \to X \}
\end{equation*}
such that for every $(x,y)\in D$ and every $t\in[0,1]\cap \Q$ it holds
\begin{equation*}
G_n(x,y)(t) \to G(x,y)(t).
\end{equation*}
Now for every $(x,y)\in D$ the function $G(x,y)$ is a $\di(x,y)$-Lipschitz function on $[0,1]\cap \Q$, in fact for every $s,t\in [0,1]\cap \Q$ it holds 
\begin{align*}
    \di \big(G(x,y)(s),G(x,y)(t)\big) &\leq \di \big(G_n(x,y)(s), G(x,y)(s)\big)+ \di \big(G_n(x,y)(s),G_n(x,y)(t)\big) \\
    &\quad + \di \big(G_n(x,y)(t), G(x,y)(t)\big)\\
    & \leq \di \big(G_n(x,y)(s), G(x,y)(s)\big)+ \di_n \big(G_n(x,y)(s),G_n(x,y)(t)\big) + \varepsilon_n \\
    &\quad + \di \big(G_n(x,y)(t), G(x,y)(t)\big)\\
    & = \di \big(G_n(x,y)(s), G(x,y)(s)\big)+ \di_n(x,y) \cdot |t-s|+ \varepsilon_n \\
    &\quad + \di \big(G_n(x,y)(t), G(x,y)(t)\big)\\
    &\to \di(x,y) \cdot |t-s|.
\end{align*}
This allows to extend $G(x,y)$ to a $\di(x,y)$-Lipschitz function on the whole interval $[0,1]$, moreover, since clearly $G(x,y)(0)=x$ and $G(x,y)(1)$, I can infer that $G(x,y)\in \Geo(X)$. Then for every $(x,y)\in D$ it is possible to extend the pointwise convergence of $G_n(x,y)$ to $G(x,y)$ to the interval $[0,1]$. In fact, given $\varepsilon>0$, for every $t \in [0,1]$ there exists $s \in [0,1]\cap \Q$ with $|t-s|<\frac{\varepsilon}{2 \di (x,y)}$ that allows to perform the following estimate
\begin{align*}
   \di \big( G_n(x,y)(t),G(x,y)(t) \big) &\leq \di \big( G_n(x,y)(t),G_n(x,y)(s) \big) + \di \big( G_n(x,y)(s),G(x,y)(s) \big) \\
   & \quad + \di \big( G(x,y)(s),G(x,y)(t) \big) \\
   &\leq (\di_n (x,y)+\di(x,y)) \cdot |t-s| + \di \big( G_n(x,y)(s),G(x,y)(s) \big)\\
   & \to 2 \di(x,y) \cdot |t-s| < \varepsilon,
\end{align*}
the claim follows from the arbitrariness of $\varepsilon$. I end up with the map
\begin{equation*}
    G: (D, \di_2) \to \big( \Geo (X), \norm{\cdot}_\text{sup} \big) \subset \big(C([0,1],X), \norm{\cdot}_\text{sup} \big)
\end{equation*}
that is the pointwise limit of the $L$-Lipschitz maps $G_n$, thus also $G$ is an $L$-Lipschitz map. Therefore it can be extended to an $L$-Lipschitz function on the whole space $X \times X$, furthermore, since $\Geo(X)$ is closed with respect to the sup norm, I obtain
\begin{equation*}
    G: (X \times X, \di_2) \to \big( \Geo (X), \norm{\cdot}_\text{sup} \big).
\end{equation*}
Then the equicontinuity of the maps $G_n$ ensures that the sequence $\{G_n\}_{n\in \setN}$ uniformly converges to $G$. 

In order to conclude the proof I only need to show that $G$ is a consistent geodesic flow, proving property 3 of Definition \ref{def:nicegeodesicflow}. To this aim fix $(x,y)\in X\times X$, $s,t\in[0,1]$ and a small $\varepsilon>0$, subsequently take $n\in \setN$ such that $\norm{G_n-G}_\text{sup}<\varepsilon$. Then it holds that
\begin{align*}
    \big\lVert \rest{s}{t} G(x,y) &- G\big(G(x,y)(s),G(x,y)(t)\big) \big\rVert & \\
    &\leq 2 \norm{G_n-G}_\text{sup} + \norm{\rest{s}{t} G_n(x,y) - G_n\big(G(x,y)(s),G(x,y)(t)\big)} \\
    &=  2 \norm{G_n-G}_\text{sup} + \norm{G_n\big(G_n(x,y)(s),G_n(x,y)(t)\big) - G_n\big(G(x,y)(s),G(x,y)(t)\big)}\\
    &\leq (2+2L) \cdot \norm{G_n-G}_\text{sup}< (2+2L)\varepsilon,
\end{align*}
thesis follows from the arbitrariness of $\varepsilon$.
\end{pr}

\noindent Once Proposition \ref{prop:ngf} has provided a consistent geodesic flow for the limit space, the next result shows how, under suitable assumptions, it is possible to prove the very strict CD condition for the metric measure space $(X,\di,\m)$.

\begin{prop}\label{prop:stabilitycompact}
Under the same assumptions of the last proposition, suppose that there exists a consistent plan selection $\Pi$ on $(X,\di,\m)$, associated to $G$, such that for every $\mu,\nu\in\Prob_{ac}(X)$ there exists a sequence $\pi_n\in \OptPlans_{\di_n}(\mu,\nu)$ satisfying
\begin{enumerate}
    \item $\pi_n \weakto \Pi(\mu,\nu)$ (up to the extraction of a subsequence),
    \item the $K$-convexity of the entropy functional holds along the $\di_n$-Wasserstein geodesic $(G_n)_\# \pi_n$, with respect to the distance $\di_n$.
\end{enumerate}
Then the metric measure space $(X,\di,\m)$ is a very strict $\CD(K,\infty)$ space.
\end{prop}

\begin{pr}
 Fix a time $t \in [0,1]$ and notice that the assumption 2 ensures that 
 \begin{equation}\label{eq:entconvG_n}
     \Ent\big( \big[(G_n)_t\big]_\# \pi_n\big) \leq (1-t) \Ent(\mu) + t \Ent(\nu) - \frac K2 t(1-t) (W_2^{\di_n})^2(\mu,\nu).
 \end{equation}
 Now, since in compact space weak convergence and Wasserstein convergence coincide, it holds that $W_2^2(\pi_n,\Pi(\mu,\nu))\to 0$. Then taking an optimal transport plan $\eta$ between $\pi_n$ and $\Pi(\mu,\nu)$ and having in mind that $G$ is $L$-Lipschitz, it is possible to do the following estimate
 \begin{align*}
     W^2_2 \big( (G_t)_\# \pi_n, (G_t)_\# \Pi(\mu,\nu)\big) &\leq \int \di^2 \big(G_t(x_1,y_1), G_t(x_2,y_2)\big) \de \eta \big((x_1,y_1),(x_2,y_2)\big)\\
     & \leq \int L^2 \cdot \di_2^2 \big((x_1,y_1),(x_2,y_2)\big)\de \eta \big((x_1,y_1),(x_2,y_2)\big) \\
     &= L^2 \cdot W_2^2(\pi_n,\Pi(\mu,\nu))\to 0.
 \end{align*}
 Consequently, I am able to infer that
 \begin{align*}
     W^2_2 \big( \big[(G_n)_t\big]_\# \pi_n, (G_t)_\# \Pi(\mu,\nu)\big) &\leq 2 W^2_2 \big( \big[(G_n)_t\big]_\# \pi_n, (G_t)_\# \pi_n\big)+ 2 W^2_2 \big( (G_t)_\# \pi_n, (G_t)_\# \Pi(\mu,\nu)\big)\\
     &\leq 2 \int \di^2 \big((G_n)_t, G_t\big) \de \pi_n + 2 W^2_2 \big( (G_t)_\# \pi_n, (G_t)_\# \Pi(\mu,\nu)\big)\\
     &\leq 2 \norm{G_n-G}_\text{sup}^2 + + 2 W^2_2 \big( (G_t)_\# \pi_n, (G_t)_\# \Pi(\mu,\nu)\big) \to 0,
 \end{align*}
 and thus that $\big[(G_n)_t\big]_\# \pi_n \xrightarrow{W_2} (G_t)_\# \Pi(\mu,\nu)$. Finally, since obviously $(W_2^{\di_n})^2(\mu,\nu) \to W_2^2(\mu,\nu)$, it is possible to pass to the limit in \eqref{eq:entconvG_n} using the lower semicontinuity of the entropy and obtain 
 \begin{equation*}
     \Ent\big( (G_t)_\# \Pi(\mu,\nu)\big) \leq (1-t) \Ent(\mu) + t \Ent(\nu) - \frac K2 t(1-t) W^2_2(\mu,\nu),
 \end{equation*}
 which, associated to Theorem \ref{thm:verystrict}, allows to conclude the proof because $t$ is arbitrary.
\end{pr}



Following verbatim the proof of Proposition \ref{prop:stabilitycompact} it is easy to deduce the following slight generalization. 

\begin{corollary}
Under the same assumptions of Proposition \ref{prop:ngf}, suppose that there exists a consistent plan selection $\Pi$ on $(X,\di,\m)$, associated to $G$. Moreover assume that for every $\mu,\nu\in\Prob_{ac}(X)$ there exist three sequences $\mu_n,\nu_n \in \Prob_{ac}(X)$ and $\pi_n\in \OptPlans_{\di_n}(\mu_n,\nu_n)$ satisfying
\begin{enumerate}
    \item $\mu_n \weakto \mu$, $\nu_n \weakto \nu$ and $\Ent(\mu_n)\to \Ent(\mu)$, $\Ent(\nu_n)\to\Ent(\nu)$,
    \item $\pi_n \weakto \Pi(\mu,\nu)$ (up to the extraction of a subsequence),
    \item the $K$-convexity of the entropy functional holds along the $\di_n$-Wasserstein geodesic $(G_n)_\# \pi_n$.
\end{enumerate}
Then the metric measure space $(X,\di,\m)$ is a very strict $\CD(K,\infty)$ space.
\end{corollary}


As already anticipated before, similar results can be proven for suitable measured Gromov Hausdorff limit spaces, also in the non-compact case. These generalizations require some technical assumption but their proof basically follow the proofs I have just presented. Anyway, in order to be concise, I prefer not to present the most general statements, except for the following proposition, which will be fundamental in the next section. The reader can easily notice that it can be proven following the proof of Proposition \ref{prop:stabilitycompact}, except for two technical details that I will fix below.

\begin{prop}\label{prop:stabilitylast}
Let $(X,\di, \m)$ be a locally compact metric measure space and let $\{\di_n\}_{n\in \setN}$ be a sequence of distances on $X$ (inducing the same topology), locally uniformly convergent to $\di$ as $n\to \infty$, such that there exists a constant $H$ satisfying 
\begin{equation}\label{eq:weirdcondition}
    \di_n(x,y) \leq H \di(x,y) \qquad \text{for every }x,y \in X\times X,
\end{equation}
for every $n$. Assume that there exists a map $G:X \times X \to C([0,1],X)$ which is a Lipschitz consistent geodesic flow for $\di$ and a consistent geodesic flow for every distance $\di_n$. 
Moreover, suppose that there is a consistent plan selection $\Pi$ on $(X,\di,\m)$, associated to $G$, such that for every $\mu,\nu\in\Prob_{ac}(X)$ there exists a sequence $\pi_n\in \OptPlans_{\di_n}(\mu,\nu)$, satisfying
\begin{enumerate}
    \item $\pi_n \weakto \Pi(\mu,\nu)$ (up to the extraction of a subsequence),
    \item the $K$-convexity of the entropy functional holds along the $\di_n$-Wasserstein geodesic $G_\# \pi_n$, with respect to the distance $\di_n$.
\end{enumerate}
Then the metric measure space $(X,\di,\m)$ is a very strict $\CD(K,\infty)$ space.
\end{prop}

\begin{remark}
Notice that condition \eqref{eq:weirdcondition} ensures that $\Prob_{ac}(\R^N,\di)\subseteq \Prob_{ac}(\R^N,\di_n)$ for every $n$.
\end{remark}

\begin{pr}
In order to repeat the same strategy used for Proposition \ref{prop:stabilitycompact} I only need to prove that $W_2^2(\pi_n,\Pi(\mu,\nu))\to 0$ and that $\lim_{n \to \infty}(W_2^{\di_n})^2(\mu,\nu) = W_2^2(\mu,\nu)$. For the first condition, according to Proposition \ref{prop:chrW2covergence}, it is sufficient to prove that \begin{equation*}\int \di_2^2 \big((x,y),(x_0,y_0)\big) \de \pi_n (x,y) \to \int \di_2^2 \big((x,y),(x_0,y_0)\big) \de \Pi(\mu,\nu) (x,y),\end{equation*}for every fixed $(x_0,y_0)\in X\times X$. But this can be easily shown, in fact for every $n\in \setN$ it holds\begin{align*}    \int \di_2^2 \big((x,y),(x_0,y_0)\big) \de \pi_n (x,y) &= \int \big[ \di^2(x,x_0) + \di^2(y,y_0) \big] \de \pi_n (x,y) \\ &= \int \di^2(x,x_0) \de \mu (x) + \int \di^2(y,y_0) \de \nu (y)\\ &= \int \di_2^2 \big((x,y),(x_0,y_0)\big) \de \Pi(\mu,\nu) (x,y).\end{align*}
On the other hand, taking $\pi\in \OptPlans_\di(\mu,\nu)$, condition \eqref{eq:weirdcondition} allows to use the dominated convergence theorem and deduce
\begin{equation*}
    \limsup_{n \to \infty}(W_2^{\di_n})^2(\mu,\nu) \leq \limsup_{n \to \infty} \int \di_n^2 \de \pi = \int \di^2 \de \pi  = W_2^2(\mu,\nu).
\end{equation*}
Moreover for every compact set $K\subset X\times X$ there exists a continuous function $\phi_K:X\times X\to [0,1]$ such that $\phi_K=0$ outside a compact set $K'$ and $f_K\equiv 1$ on $K$. Then $\phi_K \di^2_n \to \phi_K \di^2$ uniformly, therefore

\begin{equation*}
    \liminf_{n \to \infty}(W_2^{\di_n})^2(\mu,\nu) \geq \liminf_{n \to \infty}\int \phi_K \di_n^2 \de \pi_n = \int \phi_K \di^2 \de \pi \geq \int_K \di^2 \de \pi.
\end{equation*}
Since $K$ is arbitrary it is possible to conclude that 
\begin{equation*}
    \liminf_{n \to \infty}(W_2^{\di_n})^2(\mu,\nu) \geq  \int \di^2 \de \pi = W_2^2(\mu,\nu),
\end{equation*}
and consequently that $\lim_{n \to \infty}(W_2^{\di_n})^2(\mu,\nu) = W_2^2(\mu,\nu)$.
Having that $W_2^2(\pi_n,\Pi(\mu,\nu))\to 0$ and that $\lim_{n \to \infty}(W_2^{\di_n})^2(\mu,\nu) = W_2^2(\mu,\nu)$, the proof of Proposition \ref{prop:stabilitycompact} can be repeated step by step and gives the thesis.
\end{pr}

\begin{remark}
This section has shown how the existence of a consistent geodesic flow and a consistent plan selection associated to it, can help in proving the very strict CD condition. However, I have not stated any results (except for Proposition \ref{prop:ngf}) that would guarantee the existence of these two objects in a metric measure space. To this aim, it would be very interesting to investigate under which assumptions on a given consistent geodesic flow $G$ (or on the metric measure space), there exists a consistent plan selection associated to $G$. In the next section I will show how a (double) minimization procedure allows to identify a consistent plan selection in a particular metric measure space. It is possible that these arguments can also apply to a more general context.
\end{remark}

\section{Application to Crystalline Norms in $\R^N$}

The aim of this section is to prove the very strict $\CD(0,\infty)$ condition for $\R^N$ equipped with a crystalline norm and with the Lebesgue measure, using the theory developed in the last section and in particular Proposition \ref{prop:stabilitylast}. Let me point out that the Optimal Transport problem in these particular metric spaces has been already studied by Ambrosio Kirchheim and Pratelli in \cite{ambkirpra03}. They were able to solve the $L^1$-Monge problem using a secondary variational minimization in order to suitably decompose the space in transport rays. Despite the problem I want to face and the way I will do it are different from the theory developed in \cite{ambkirpra03}, I will in turn use a secondary variational problem to select a suitable transport plan connecting two given measures, obtaining, as a byproduct, the existence of optimal transport map between them.

Before going on, I fix the notation I will use in this section.
Given a finite set of vectors $\tilde {\mathcal V}\subset \R^N$ such that $\text{span}(\tilde{\mathcal V})=\R^N$, introduce the associate crystalline norm, which is defined as follows
\begin{equation*}
    \norm{x}:= \max_{v\in \tilde{\mathcal{V}}} |\scal{x}{v}|
\end{equation*}
and the corresponding distance 
\begin{equation*}
    \di(x,y):= \norm{x-y}= \max_{v\in \tilde{\mathcal{V}}}|\scal{x-y}{v}|.
\end{equation*}
For sake of exposition, from now on I am going to use the following equivalent formulations for the norm and the distance:
\begin{equation*}
    \norm{x}:= \max_{v\in \mathcal{V}} \scal{x}{v}, \quad \di(x,y):= \norm{x-y}= \max_{v\in \mathcal{V}}\scal{x-y}{v},
\end{equation*}
where $\mathcal V$ denotes the set $\tilde{\mathcal V} \cup (-\tilde{\mathcal V})$.

As the reader can easily guess, in this framework the choice of a consistent geodesic flow is not really problematic, in fact it is sufficient to consider the Euclidean one, that is
\begin{equation*}
\begin{split}
      G: \R^N \times \R^N &\to C([0,1],\R^N) \\
    (x,y) &\mapsto (t \mapsto (1-t)x + t y)
\end{split}.
\end{equation*}
The rest of the chapter will be then dedicated to the choice of a suitable plan selection, associated to $G$, satisfying the requirements of Proposition \ref{prop:stabilitylast}. It will be identified via a secondary variational minimization. This type of procedure turns out to be useful in many situation (see for example Chapter 2 and 3 in \cite{santambrogio2015optimal}) and in this specific case is inspired by the work of Rajala \cite{rajala2013failure}. Let me now go into the details.
Given two measures $\mu,\nu\in\ProbTwo(\R^N)$, consider the usual Kantorovich problem with cost $c(x,y)=\di^2(x,y)$, that is 
\begin{equation*}
    \min_{\pi\in\Gamma(\mu,\nu)} \int_{\R^N\times \R^N} \di^2(x,y) \de \pi(x,y),
\end{equation*}
calling $\Pi_1(\mu,\nu)$ the set of its minimizers. Consequently consider the secondary variational problem 
\begin{equation}\label{eq:2ndvarproblem}
    \min_{\pi\in\Pi_1(\mu,\nu)} \int_{\R^N\times \R^N} \di^2_{eu}(x,y) \de \pi(x,y),
\end{equation}
where I denote by $\di_{eu}$ the Euclidean distance, and denote by $\Pi_2(\mu,\nu)\subseteq\Pi_1(\mu,\nu)$ the set of minimizers, which can be easily seen to be not empty. In Theorem \ref{thm:crystalline} I will show that, if $\mu$ is absolutely continuous, $\Pi_2(\mu,\nu)$ consists of a single element, but, in order to do this I have to preliminarily exploit the cyclical monotonicity properties of the plans in $\Pi_2(\mu,\nu)$.

\begin{prop}\label{prop:doublemonotonicity}
Every $\pi\in \Pi_2(\mu,\nu)$ is concentrated in a set $\Gamma$, such that for every $(x,y),(x',y')\in \Gamma$ it holds that
\begin{equation}\label{eq:monot1}
    \di^2(x,y)+\di^2(x',y') \leq \di^2(x,y')+ \di^2(x',y),
\end{equation}
moreover, if $\di^2(x,y)+\di^2(x',y') = \di^2 (x,y')+ \di^2(x',y)$, then 
\begin{equation}\label{eq:monot2}
    \di^2_{eu}(x,y)+\di^2_{eu}(x',y') \leq \di^2_{eu} (x,y')+ \di^2_{eu}(x',y).
\end{equation}
\end{prop}

\begin{pr}
Fix $\pi\in \Pi_2(\mu,\nu)$ and notice that, since in particular $\pi\in \Pi_1(\mu,\nu)$, Proposition \ref{prop:cmonotonicity} yields that $\pi$ is concentrated in a set $\Gamma_1$ satisfying \eqref{eq:monot1}. Furthermore, according to Proposition \ref{prop:Kantorovichpotential} and Remark \ref{remark:Kantorovich}, fix an upper semicontinuous Kantorovich potential $\phi$ for the cost $c(x,y)=\di^2(x,y)$, such that also $\phi^c$ is upper semicontinuous. In particular for every $\eta\in \Pi_1(\mu,\nu)$, it holds 
\begin{equation*}
    \phi(x)+\phi^c(y)= c(x,y)=\di^2(x,y), \quad \text{for $\eta$-almost every $(x,y)\in \R^N\times\R^N$}.
\end{equation*}
As a consequence, notice that being a minimizer of the secondary variational problem \eqref{eq:2ndvarproblem} is equivalent to realize the minimum of 
\begin{equation*}
     \min_{\eta\in\Pi(\mu,\nu)} \int_{\R^N\times \R^N} \tilde c(x,y) \de \eta(x,y),
\end{equation*}
where the cost $\tilde c$ is defined as
\begin{equation*}
\tilde c(x,y)=
    \begin{cases}
  \di_{eu}^2(x,y) &\text{if } \phi(x)+ \phi^c(y)= \di^2(x,y) \\
       +\infty &\text{otherwise}
\end{cases}.
\end{equation*}
Observe that, since $\phi$ and $\phi^c$ are upper semicontinuous, the cost $\tilde c$ is lower semicontinuous. Thus Proposition \ref{prop:cmonotonicity} ensures that $\pi$ is concentrated in a set $\Gamma_2$ which is $\tilde c$-cyclically monotone. Moreover, up to modify $\Gamma_2$ in a $\pi$-null set, it is possible to assume that for every $(x,y)\in\Gamma_2$
\begin{equation*}
     \phi(x)+\phi^c(y)= c(x,y)=\di^2(x,y).
\end{equation*}
Now take $(x,y),(x',y')\in \Gamma_2$ with $\di^2(x,y)+\di^2(x',y') = \di^2 (x,y')+ \di^2(x',y)$ and deduce that 
\begin{align*}
    \phi(x)+\phi^c(y)+\phi(x')+\phi^c(y')=\di^2(x,y)+\di^2(x',y') = \di^2 (x,y')+ \di^2(x',y).
\end{align*}
On the other hand $\phi(x)+\phi^c(y')\leq \di^2 (x,y')$ and $\phi(x')+\phi^c(y)\leq \di^2 (x',y)$, therefore I obtain 
\begin{equation*}
    \phi(x)+\phi^c(y')= \di^2 (x,y') \quad \text{and}\quad \phi(x')+\phi^c(y)= \di^2 (x',y).
\end{equation*}
Finally the $\tilde c$-cyclical monotonicity allows to conclude that
\begin{align*}
     \di^2_{eu}(x,y)+\di^2_{eu}(x',y')&= \tilde c(x,y)+\tilde c(x',y')\\
     &\leq  \tilde c(x,y')+\tilde c(x',y)= \di^2_{eu} (x,y')+ \di^2_{eu}(x',y),
\end{align*}
which is exactly \eqref{eq:monot2}. Summing up, it is easy to check that the set $\Gamma=\Gamma_1 \cap \Gamma_2$ satisfies the requirements of Proposition \ref{prop:doublemonotonicity}.
\end{pr}

\noindent I can now go into the proof of one of the main results of this work.

\begin{theorem}\label{thm:crystalline}
Given two measures $\mu,\nu\in \ProbTwo(\R^N)$ with $\mu$ absolutely continuous with respect to $\Leb^n$, there exists a unique $\pi\in \Pi_2(\mu,\nu)$ and it is induced by a map. 
\end{theorem}

\begin{pr}
Reasoning as in Remark \ref{rem:uniquenessthroughmap}, it is sufficient to prove that every plan in $\Pi_2(\mu,\nu)$ is induced by a map.
So take $\pi\in \Pi_2(\mu,\nu)$, applying Proposition \ref{prop:doublemonotonicity} it is possible to find a full $\pi$-measure set $\Gamma$, satisfying the monotonicity requirements \eqref{eq:monot1} and \eqref{eq:monot2}. Assume by contradiction that $\pi$ is not induced by a map, calling $(\pi_x)_{x\in \R^N}\subset \Prob(\R^N)$ the disintegration with respect to the projection map $\p_1$, then $\pi_x$ is not a delta measure for a $\mu$-positive set. Moreover, given a non-empty set $V \subseteq \mathcal V$, define the sets
\begin{equation*}
     \tilde A_{z,V}:= \big\{x\in \R^N \suchthat \di(z,x)=\scal{z-x}{v}\,\,\, \text{for every } v \in V\big\},
\end{equation*}
\begin{equation*}
    A_{z,V} := \big\{x \in\tilde A_{z,V} \suchthat \di(z,x)>\scal{z-x}{v}\,\,\, \text{for every } v \in\mathcal V \setminus V \big\},
\end{equation*}
\begin{equation*}
    A_{z,V}^\varepsilon := \big\{x \in\tilde A_{z,V} \suchthat \di(z,x)>\scal{z-x}{v} + \varepsilon\,\,\, \text{for every } v \in\mathcal V \setminus V \big\}.
\end{equation*}
Notice that, for every fixed $z\in\R^N$, the sets $A_{z,V}$ constitute a partition of $\R^N$ as $V\subseteq \mathcal V$ varies.
Consequently, I divide the proof in three steps, whose combination will allow me to conclude by contradiction.

\noindent \textit{Step 1: Given two nonempty sets $V_1, V_2 \subseteq \mathcal V$ such that $v_1\ne v_2$ for every $v_1\in V_1$ and $v_2\in V_2$ (that is $V_1 \cap V_2=\emptyset$), the set
\begin{equation*}
    E:= \big\{z\in \R^N \suchthat \pi_z(A_{z,V_1})>0 \text{ and } \pi_z(A_{z,V_2})>0 \big\}
\end{equation*}
has zero $\mu$-measure.}

First of all, notice that if $E$ is non-empty, then for every fixed $z\in \R^N$ there exist $x\in A_{z,V_1}$ and $y\in A_{z,V_2}$ such that $\di(z,x)=\di(z,y)=1$ and in particular  
\begin{equation*}
    \scal{x}{v_1}=1 > \scal{y}{v_1} \quad \text{for every }v_1\in V_1
\end{equation*}
and
\begin{equation*}
    \scal{y}{v_2}=1 > \scal{x}{v_2} \quad \text{for every }v_2\in V_2.
\end{equation*}
Therefore, calling $\bar v =x-y$, it holds that
\begin{equation}\label{eq:barv}
    \begin{split}
        &\scal{\bar v}{v_1}>0 \quad \text{for every }v_1\in V_1,\\
        &\scal{\bar v}{v_2}<0 \quad \text{for every }v_1\in V_2.
    \end{split}
\end{equation}

Now, assume by contradiction that $E$ has positive $\mu$-measure, in particular it is non-empty and there exists $\bar v$ satisfying \eqref{eq:barv}. Moreover, notice that, since $\Gamma$ is $\pi$-measurable and has full measure, then $\Gamma_z:=\{(z,y)\in \Gamma : y\in \R^N\}$ is $\pi_z$-measurable with $\pi_z(\Gamma_z)=1$ for $\mu$-almost every $z\in \R^N$. In particular for $\varepsilon>0$ small enough the set
\begin{equation*}
    E_\varepsilon:= \big\{z\in \R^N \suchthat \pi_z(A_{z,V_1}^\varepsilon \cap \Gamma_z)>0 \text{ and } \pi_z(A_{z,V_2}^\varepsilon \cap \Gamma_z)>0 \big\}
\end{equation*}
has positive $\mu$-measure, and thus it also has positive $\Leb^N$-measure. Take a Lebesgue density point $\bar z$ of $E_\varepsilon$, then in a neighborhood of $\bar z$ there exist $z$ such that $z,z+\epsilon\bar v \in E_\varepsilon$ for a suitable $0<\epsilon<\frac{ \varepsilon}{ \norm{\bar v}}$. Now, there exist $x\in A_{z,V_1}^\varepsilon$ and $y\in A_{z+\epsilon \bar v,V_2}^\varepsilon$ such that $(z,x),(z+\epsilon\bar v,y) \in \Gamma$. Notice that for every $v_1\in V_1$, it holds
\begin{equation}\label{eq:distance1}
    \scal{x - (z + \epsilon \bar v)}{v_1} = \scal{x-z}{v_1} -\epsilon \scal{\bar v}{v_1} < \scal{x-z}{v_1} = \di(z,x) ,
\end{equation}
while for every $w\in \mathcal V \setminus  V_1$ it is possible to perform the following estimate:
\begin{equation}\label{eq:distance2}
    \scal{x - (z + \epsilon \bar v)}{w} = \scal{x-z}{w} -\epsilon \scal{\bar v}{w}< \di (x,z) - \varepsilon + \epsilon \norm{\bar v} < \di (z,x).
\end{equation}
The combination of \eqref{eq:distance1} and \eqref{eq:distance2} yields 
\begin{equation}\label{eq:ineq1}
    \di(x,z+\epsilon \bar v) < \di (z,x).
\end{equation}
Similarly, it holds 
\begin{equation*}
    \scal{y-z}{v_2} = \scal{y-(z+\epsilon\bar v)}{v_2} +\epsilon \scal{\bar v}{v_2} < \scal{y-(z+\epsilon\bar v)}{v_2} = \di (z + \epsilon \bar v,y) ,
\end{equation*}
for every $v_2\in V_2$, and 
\begin{equation*}
    \scal{y-z}{w} = \scal{y- (z+\epsilon \bar v)}{w} + \epsilon \scal{\bar v}{w}< \di (z + \epsilon \bar v,y) - \varepsilon + \epsilon \norm{\bar v} < \di (z + \epsilon \bar v,y),
\end{equation*}
for every $w \in \mathcal V \setminus V_2$, which together show that 
\begin{equation}\label{eq:ineq2}
      \di (z,y)<\di(z+\epsilon \bar v,y).
\end{equation}
Now, the inequalities \eqref{eq:ineq1} and \eqref{eq:ineq2} allow to infer that
\begin{equation*}
    \di^2 (z,x) + \di^2(z+\epsilon\bar v, y) > \di^2 (z,y) + \di^2 (z+\epsilon\bar v, x),
\end{equation*}
contradicting the condition \eqref{eq:monot1} of Proposition \ref{prop:doublemonotonicity}.

\noindent \textit{Step 2: Given two nonempty sets $V_1, V_2 \subseteq \mathcal V$ such that $V_1 \cap V_2 \ne \emptyset$ and $V_1 \ne V_2$, the set
\begin{equation*}
    E:= \big\{z\in \R^N \suchthat \pi_z(A_{z,V_1})>0 \text{ and } \pi_z(A_{z,V_2})>0 \big\}
\end{equation*}
has zero $\mu$-measure.}\\
 Call $V=V_1 \cap V_2$, $W_1=V_1 \setminus V$ and $W_2=V_2 \setminus V$. Assume by contradiction that $E$ has positive $\mu$-measure, then for $\varepsilon>0$ sufficiently small the set
\begin{equation*}
    E_\varepsilon:= \big\{z\in \R^N \suchthat \pi_z(A_{z,V_1}^\varepsilon)>0 \text{ and } \pi_z(A_{z,V_2}^\varepsilon)>0 \big\}
\end{equation*}
has positive $\mu$-measure too. As a consequence 
\begin{equation*}
    \gamma:= \int_{E_\varepsilon} \restr{\pi_z}{A_{z,V_1}^\varepsilon} \times \restr{\pi_z}{A_{z,V_2}^\varepsilon} \de \mu (z)
\end{equation*}
is a strictly positive measure on $\R^N\times\R^N$ with $\gamma \big(\big\{(x,x):x\in \R^N\big\}\big)=0$. Thus there exists $(\bar x, \bar y)\in \supp(\gamma)$ with $\bar x \ne \bar y$ and then 
\begin{equation*}
    \gamma \big( B_{\delta}(\bar x)\times B_{\delta}(\bar y)\big) >0,
\end{equation*}
for every $\delta>0$. In particular, proceeding as in the first step, it is possible to conclude that for every $\delta>0$ the set  
\begin{equation*}
    E_\varepsilon^\delta:= \big\{z\in \R^N \suchthat \pi_z\big(A_{z,V_1}^\varepsilon \cap \Gamma_z \cap B_{\delta}(\bar x)\big)>0 \text{ and } \pi_z\big(A_{z,V_2}^\varepsilon \cap \Gamma_z\cap B_{\delta}(\bar y)\big)>0 \big\}
\end{equation*}
has positive $\mu$-measure, and thus it also has positive $\Leb^N$-measure. Now, I divide the proof in two cases, depending on the vector $\bar v = \bar x-\bar y$:
\begin{itemize}
    \item Case 1: $\scal{\bar v}{v}=0$ for every $v \in V$.\\
    Since $(\bar x,\bar y)\in \supp (\gamma)$, for every $\eta>0$ there exist $x_\eta,y_\eta, z_\eta$ such that $\norm{\bar x - x_\eta},\norm{\bar y - y_\eta}<\eta$ and $x_\eta \in A_{z_\eta,V_1}^\varepsilon$, $y_\eta \in A_{z_\eta,V_2}^\varepsilon$. Then, given $v\in V$, for every $v_1\in W_1$ it holds that 
    \begin{align*}
        \scal{x_\eta- z_\eta}{v_1} = \scal{x_\eta- z_\eta}{v} &= \scal{x_\eta- \bar x}{v}+ \scal{\bar v}{v}+ \scal{\bar y- y_\eta}{v}+ \scal{y_\eta- z_\eta}{v}\\
        &> \scal{y_\eta- z_\eta}{v_1} +\frac \varepsilon 2,
    \end{align*}
    for $\eta$ small enough. Thus, if $\eta$ is sufficiently small, follows that 
    \begin{equation*}
        \scal{x_\eta}{v_1} > \scal{y_\eta}{v_1} +\frac \varepsilon 2 \quad \text{for every } v_1 \in W_1,
    \end{equation*}
    and similarly 
    \begin{equation*}
        \scal{x_\eta}{v_2} < \scal{y_\eta}{v_2} - \frac \varepsilon 2 \quad \text{for every } v_2 \in W_2.
    \end{equation*}
    Taking the limit as $\eta \to 0$, clearly $x_\eta \to \bar x$ and $y_\eta \to \bar y$, therefore I conclude that 
     \begin{equation}\label{eq:V1}
        \scal{\bar x}{v_1} > \scal{\bar y}{v_1} \, \text{ and thus } \, \scal{\bar v}{v_1}>0,  \quad \text{for every } v_1 \in W_1,
    \end{equation}
    and 
    \begin{equation}\label{eq:V2}
        \scal{\bar x}{v_2} < \scal{\bar y}{v_2} \, \text{ and thus } \, \scal{\bar v}{v_2}<0, \quad \text{for every } v_2 \in W_2.
    \end{equation}
    
    Now, fix $\delta>0$ sufficiently small such that 
    \begin{equation}\label{eq:choicedelta}
        \scal{\bar v}{x}> \scal{\bar v}{y}, \quad \text{for every } x \in B_\delta (\bar x) \text{ and } y \in B_\delta (\bar y). 
    \end{equation}
    As already emphasized, the set $E_\varepsilon^\delta$ has positive Lebesgue measure, then take one of its density points $\bar z$. In a neighborhood of $\bar z$ there exists $z$, such that $z, z + \epsilon \bar v \in E_\varepsilon^\delta$ for a suitable $0<\epsilon < \frac{\varepsilon}{\norm{\bar v}}$, subsequently take $x\in A_{z,V_1}^\varepsilon\cap B_\delta(\bar x)$ with $(z,x)\in \Gamma$, and $y\in A_{z+\epsilon \bar v,V_2}^\varepsilon \cap B_\delta(\bar y)$ with $(z+\epsilon \bar v,y)\in \Gamma$. Notice that for every $v\in V$ it holds 
    \begin{equation*}
        \scal{x-(z+\epsilon \bar v)}{v} = \scal{x-z}{v}= \di (z,x), 
    \end{equation*}
    moreover \eqref{eq:V1} ensures that for every $v_1 \in W_1$ 
    \begin{equation*}
        \scal{x-(z+\epsilon \bar v)}{v_1} < \scal{x-z}{v_1}= \di (z,x), 
    \end{equation*}
    while for every $w\in \mathcal V \setminus V_1$ the following estimate can be performed 
    \begin{equation*}
         \scal{x-(z+\epsilon \bar v)}{w} = \scal{x-z}{w}- \epsilon \scal{\bar v}{w}<  \di (x,z) -\varepsilon + \epsilon \norm{\bar v}  < \di (z,x). 
    \end{equation*}
    This last three relations show that 
    \begin{equation}\label{eq:equalitycrystalline1}
        \di(z+\epsilon \bar v, x) = \di (z, x),
    \end{equation}
    and analogously using \eqref{eq:V2} it can be proven that 
    \begin{equation}\label{eq:equalitycrystalline2}
         \di (z, y)=\di(z+\epsilon \bar v, y).
    \end{equation}
    On the other hand, the choice of $\delta$ I made (see \eqref{eq:choicedelta}) guarantees that 
    \begin{align*}
        \di_{eu}^2(z+\epsilon \bar v, x) + \di_{eu}^2(z,y) &= \scal{z+\epsilon \bar v-x}{ z+\epsilon \bar v -x} + \scal{z-y}{z-y} \\
        &= \scal{z-x}{z-x}+2\scal{z-x}{\epsilon \bar v}+\scal{\epsilon \bar v}{\epsilon \bar v}+ \scal{z-y}{z-y}\\
        &< \scal{z-x}{z-x}+\scal{\epsilon \bar v}{\epsilon \bar v}+2\scal{z-y}{\epsilon \bar v}+ \scal{z-y}{z-y} \\
        &=  \scal{z-x}{z-x} + \scal{z+\epsilon \bar v-y}{ z+\epsilon \bar v -y} \\
        &= \di_{eu}^2 (z,x) + \di_{eu}^2(z+\epsilon \bar v, y),
    \end{align*}
    which, together with \eqref{eq:equalitycrystalline1} and \eqref{eq:equalitycrystalline2}, contradicts the condition \eqref{eq:monot2} of Proposition \ref{prop:doublemonotonicity}.
    \item Case 2: there exists $\bar w \in V$ such that $\scal{\bar v}{\bar w}\ne 0$.\\
    Without losing generality I can assume $\scal{\bar v}{\bar w}> 0$, then it is possible to fix a sufficiently small $\delta>0$ such that, for a suitable $\eta>0$, it holds
     \begin{equation*}
        \scal{\bar w}{x}> \scal{\bar w}{y}+\eta, \quad \text{for every } x \in B_\delta (\bar x) \text{ and } y \in B_\delta (\bar y). 
    \end{equation*}
    Fix a vector $\tilde v \in A_{z,V_1}$. Repeating the argument used in Case 1 it is possible to find a point $z\in \R^n$, such that $z,z+\epsilon\tilde v \in E_\varepsilon^\delta$ for a suitable $0<\epsilon< \max \big\{ \frac \varepsilon {2 \norm{\tilde v}}, \frac  \eta {2 \norm{\tilde v}} \big\}$. Then take $x\in A_{z,V_1}^\varepsilon\cap B_\delta(\bar x)$ and $y\in A_{z+\epsilon \tilde v,V_2}^\varepsilon \cap B_\delta(\bar y)$ with $(z,x),(z+\epsilon\tilde v)\in \Gamma$, and notice that for every $v_1 \in V_1$ it holds that
    \begin{equation*}
        \scal{x-(z+\epsilon \tilde v)}{v_1} = \scal{x-z}{v}  - \epsilon \scal{ \tilde v}{v_1}= \di (z,x) - \epsilon \norm{\tilde v}
    \end{equation*}
    while for every $w\in \mathcal V \setminus V_1$ I have
    \begin{equation*}
        \scal{x-(z+\epsilon \tilde v)}{w} = \scal{x-z}{w}  - \epsilon \scal{ \tilde v}{w}< \di (z,x) - \varepsilon + \epsilon \norm{\tilde v} < \di (z,x) - \epsilon \norm{\tilde v},
    \end{equation*}
   therefore follows that
   \begin{equation}\label{eq:exactnorm}
       \di (z+\epsilon\tilde v, x) = \di (z,x) - \epsilon \norm{\tilde v}.
   \end{equation}
   On the other hand, observe that 
   \begin{equation}\label{eq:doubleproduct}
       \di(z+\epsilon\tilde v, y)= \scal{y-(z+\epsilon\tilde v)}{\bar w} = \scal{y-z}{\bar w}-\epsilon\scal{\tilde v}{\bar w}<\di (z,x) - \eta  + \epsilon \norm{\tilde v}<\di (z,x).
   \end{equation}
   It is then possible to conclude that 
   \begin{align*}
       \di^2 (z+\epsilon\tilde v, x) + \di^2 (z,y) &\leq \big(\di (z,x) - \epsilon \norm{\tilde v} \big)^2 + \big(\di (z+\epsilon\tilde v, y)+ \epsilon \norm{\tilde v}\big)^2 \\
       & = \di^2 (z,x) + \di^2(z+\epsilon\tilde v, y) - 2 \epsilon \norm{\tilde v}\big( \di (z,x) - \di (z+\epsilon\tilde v, y) \big)\\
       &< \di^2 (z,x) + \di^2(z+\epsilon\tilde v, y),
   \end{align*}
   where I used both \eqref{eq:exactnorm} and \eqref{eq:doubleproduct}. This last inequality contradicts condition \eqref{eq:monot1} of Proposition \ref{prop:doublemonotonicity}.
\end{itemize}
 
\noindent  \textit{Step 3: Given a nonempty set $V \subseteq \mathcal V$, the set
\begin{equation*}
    E:= \big\{z\in \R^n \suchthat \restr{\pi_z}{A_{z,V}}\text{ is not a delta measure} \big\}
\end{equation*}
has zero $\mu$-measure.}\\
 The proof of this step is very similar to the one of Step 2, nevertheless I decided to present it anyway, but avoiding all the details which can be easily fixed following the proof of Step 2.
 Assume by contradiction that $E$ has positive $\mu$-measure, then for $\varepsilon>0$ sufficiently small the set
\begin{equation*}
    E_\varepsilon:= \big\{z\in \R^n \suchthat \restr{\pi_z}{A_{z,V}^\varepsilon}\text{ is not a delta measure} \big\}
\end{equation*}
has positive $\mu$-measure too. As a consequence 
\begin{equation*}
    \gamma:= \int_{E_\varepsilon} \restr{\pi_z}{A_{z,V}^\varepsilon} \times \restr{\pi_z}{A_{z,V}^\varepsilon} \de \mu (z)
\end{equation*}
is a strictly positive measure on $\R^N\times\R^N$ that is not concentrated on $\big\{(x,x):x\in \R^N\big\}$. Thus there exists $(\bar x, \bar y)\in \supp(\gamma)$ with $\bar x \ne \bar y$ and then 
\begin{equation*}
    \gamma \big( B_{\delta}(\bar x)\times B_{\delta}(\bar y)\big) >0,
\end{equation*}
for every $\delta>0$. In particular, proceeding as in the first step, it is possible to conclude that for every $\delta>0$ the set  
\begin{equation*}
    E_\varepsilon^\delta:= \big\{z\in \R^N \suchthat \pi_z\big(A_{z,V}^\varepsilon \cap \Gamma_z \cap B_{\delta}(\bar x)\big)>0 \text{ and } \pi_z\big(A_{z,V}^\varepsilon \cap \Gamma_z\cap B_{\delta}(\bar y)\big)>0 \big\}
\end{equation*}
has positive $\mu$-measure, and thus it also has positive $\Leb^N$-measure.
 Now, as I did in Step 2, I divide the proof in two cases:
\begin{itemize}
    \item Case 1: $\scal{\bar v}{v}=0$ for every $v \in V$.\\
    First of all, fix $\delta>0$ sufficiently small such that 
    \begin{equation*}
        \scal{\bar v}{x}> \scal{\bar v}{y}, \quad \text{for every } x \in B_\delta (\bar x) \text{ and } y \in B_\delta (\bar y). 
    \end{equation*}
    Proceeding as in Step 2, I can find $z\in \R^n$, such that $z, z + \epsilon \bar v \in E_\varepsilon^\delta$ for a positive, suitably small $\epsilon$. Subsequently take $x\in A_{z,V}^\varepsilon\cap B_\delta(\bar x)$ with $(z,x)\in \Gamma$, and $y\in A_{z+\epsilon \bar v,V}^\varepsilon \cap B_\delta(\bar y)$ with $(z+\epsilon \bar v,y)\in \Gamma$. Following the proof of Step 2, it is easy to realize that  
    \begin{equation}\label{eq:equalitycrystalline3}
        \di(z+\epsilon \bar v, x) = \di (z, x),
    \end{equation}
    and 
    \begin{equation}\label{eq:equalitycrystalline4}
         \di (z, y)=\di(z+\epsilon \bar v, y).
    \end{equation}
    On the other hand, the choice of $\delta$ I made guarantees that 
    \begin{equation*}
        \di_{eu}^2(z+\epsilon \bar v, x) + \di_{eu}^2(z,y) <\di_{eu}^2 (z,x) + \di_{eu}^2(z+\epsilon \bar v, y),
    \end{equation*}
    which, together with \eqref{eq:equalitycrystalline3} and \eqref{eq:equalitycrystalline4}, contradicts the condition \eqref{eq:monot2} of Proposition \ref{prop:doublemonotonicity}.
    \item Case 2: there exists $\bar w \in V$ such that $\scal{\bar v}{\bar w}\ne 0$.\\
    Without losing generality I can assume $\scal{\bar v}{\bar w}> 0$, then it is possible to fix a sufficiently small $\delta>0$ such that, for a suitable $\eta>0$, it holds that
     \begin{equation*}
        \scal{\bar w}{x}> \scal{\bar w}{y}+\eta, \quad \text{for every } x \in B_\delta (\bar x) \text{ and } y \in B_\delta (\bar y). 
    \end{equation*}
    Once fixed a vector $\tilde v \in A_{z,V}$, it is possible to find a point $z\in \R^n$, such that $z,z+\epsilon\tilde v \in E_\varepsilon^\delta$ for a positive, suitably small $\epsilon$. Then take $x\in A_{z,V_1}^\varepsilon\cap B_\delta(\bar x)$ and $y\in A_{z+\epsilon \tilde v,V_2}^\varepsilon \cap B_\delta(\bar y)$ with $(z,x),(z+\epsilon\tilde v)\in \Gamma$. Proceeding as I did in Step 2, it is easy to notice that 
   \begin{equation}\label{eq:exactnorm2}
       \di (z+\epsilon\tilde v, x) = \di (z,x) - \epsilon \norm{\tilde v}.
   \end{equation}
   and
   \begin{equation}\label{eq:doubleproduct2}
       \di(z+\epsilon\tilde v, y)<\di (z,x).
   \end{equation}
   Then, combining \eqref{eq:exactnorm2} and \eqref{eq:doubleproduct2}, I can conclude that 
   \begin{equation*}
       \di^2 (z+\epsilon\tilde v, x) + \di^2 (z,y) < \di^2 (z,x) + \di^2(z+\epsilon\tilde v, y),
   \end{equation*}
    contradicting condition \eqref{eq:monot1} of Proposition \ref{prop:doublemonotonicity}.
\end{itemize}
As anticipated before, it is easy to realize that the combination of the three steps allows to conclude the proof.
  \end{pr}
  
At this point it is clear that Theorem \ref{thm:crystalline} provides a plan selection on $\Prob_{ac}(\R^N)\times \Prob_{ac}(\R^N)$, simply imposing $\Pi(\mu,\nu)$ to be equal to the only optimal transport plan in $\Pi_2(\mu,\nu)$. The following proposition ensures that $\Pi$ is a consistent plan selection. 
  
\begin{prop}
The map $\Pi$ is a consistent plan selection, associated to $G$.
\end{prop}

 \begin{pr}
 Considering how $\Pi$ has been defined, in order to conclude the proof, is sufficient to prove conditions 2.1 and 2.2 of Lemma \ref{lem:conditions}. It is easy to realize that condition 2.1 is satisfied since $f \cdot \Pi(\mu,\nu)\ll \Pi(\mu,\nu)$ with bounded density, for every suitable $f$. Condition 2.2 is a little bit trickier and I am going to prove it with full details.
 
 Assume by contradiction that, for some $\mu,\nu \in \Prob_{ac}(\R^N)$, $\pi_2:= (G_s,G_t)_\# \Pi(\mu,\nu)$ is not a minimizer for the secondary variational problem \eqref{eq:2ndvarproblem}, with absolutely continuous marginals $\mu_s:=(G_s)_\# \Pi(\mu,\nu)$ and $\mu_t:=(G_t)_\# \Pi(\mu,\nu)$. Since $\pi_2$ is clearly an optimal transport plan, this means that there exists $\pi \in \OptPlans(\mu_s,\mu_t)$ such that 
 \begin{equation*}
     \int \di_{eu}^2 (x,y) \de \pi < \int \di_{eu}^2 (x,y) \de \pi_2.
 \end{equation*}
 Then Dudley's gluing lemma ensures the existence of a probability measure $\tilde \pi \in \Prob((\R^N)^4)$ such that
 \begin{equation*}
     (\p_1,\p_2)_\# \tilde \pi =\pi_1, \quad
     (\p_2,\p_3)_\# \tilde \pi = \pi \quad \text{and} \quad  (\p_3,\p_4)_\# \tilde \pi =\pi_3, 
 \end{equation*}
 where $\pi_1 := (G_0,G_s)_\# \Pi(\mu,\nu)$ and $\pi_3 := (G_t,G_1)_\# \Pi(\mu,\nu)$. Defining $\bar \pi:= (\p_1,\p_4)_\# \tilde \pi$ it is possible to perform the following estimate 
 \begin{align*}
     \int \di^2(x,y) \de \bar \pi (x,y) &= \int \di^2(x,y) \de \tilde \pi (x,z,w,y)\\
     &\leq \int \big(\di(x,z)+\di(z,w)+\di(w,y)\big)^2 \de \tilde \pi (x,z,w,y)\\
     &= \int \di^2 (x,z)  \de \pi_1 + \int \di^2 (z,w)  \de \pi + \int \di^2 (w,y)  \de \pi_3 \\
     &\text{  } + 2 \int \di(x,z) \di(z,w)\de \tilde \pi (x,z,w,y) + 2 \int \di(x,z) \di(w,y)\de \tilde \pi (x,z,w,y) \\
     & \text{  } + 2 \int \di(z,w) \di(w,y)\de \tilde \pi (x,z,w,y).
 \end{align*}
 Moreover, this last three integrals can be further estimated, inferring that
 \begin{align*}
     2 \int \di(x,z) \di(z,w)\de \tilde \pi (x,z,w,y) &= s (t-s) \int 2 \bigg(\frac{1}{s} \di(x,z) \bigg) \bigg( \frac{1}{t-s}\di(z,w) \bigg) \de \tilde \pi (x,z,w,y) \\ 
     &\leq \frac{t-s}{s} \int \di^2 (x,z)  \de \pi_1 + \frac{s}{t-s} \int \di^2 (z,w)  \de \pi 
 \end{align*}
 and similarly 
 \begin{equation*}
     2 \int \di(x,z) \di(w,y)\de \tilde \pi (x,z,w,y) \leq \frac{1-t}{s} \int \di^2 (x,z)  \de \pi_1 + \frac{s}{1-t} \int \di^2 (w,y)  \de \pi,
 \end{equation*}
 \begin{equation*}
     2 \int \di(z,w) \di(w,y)\de \tilde \pi (x,z,w,y) \leq \frac{1-t}{t-s} \int \di^2 (z,w)  \de \pi_1 + \frac{t-s}{1-t} \int \di^2 (w,y)  \de \pi.
 \end{equation*}
 Putting together this last three inequalities, it is possible to deduce that 
 \begin{align*}
     \int \di^2(x,y) \de \bar \pi (x,y) &\leq \frac{1}{s}  \int \di^2 (x,z)  \de \pi_1 + \frac{1}{t-s}  \int \di^2 (z,w)  \de \pi + \frac{1}{1-t}  \int \di^2 (w,y)  \de \pi_3 \\
     &= \frac 1s W_2^2(\mu,\mu_s) + \frac{1}{t-s} W_2^2(\mu_s,\mu_t) + \frac{1}{1-t} W_2^2 (\mu_t,\nu) = W_2^2(\mu,\nu),
 \end{align*}
 where I used the fact that $G_\# \Pi(\mu,\nu)$ is an optimal geodesic plan. In particular this shows that $\bar \pi \in \OptPlans(\mu,\nu)$. Furthermore, performing the same computation as before, one can infer that 
  \begin{align*}
     \int \di_{eu}^2(x,y) \de \bar \pi (x,y) &\leq \frac{1}{s}  \int \di_{eu}^2 (x,z)  \de \pi_1 + \frac{1}{t-s}  \int \di_{eu}^2 (z,w)  \de \pi + \frac{1}{1-t}  \int \di_{eu}^2 (w,y)  \de \pi_3 \\
     &< \frac{1}{s}  \int \di_{eu}^2 (x,z)  \de \pi_1 + \frac{1}{t-s}  \int \di_{eu}^2 (z,w)  \de \pi_2 + \frac{1}{1-t}  \int \di_{eu}^2 (w,y)  \de \pi_3 \\
     &= \int \di_{eu}^2(x,y) \de \Pi(\mu,\nu),
 \end{align*}
 where this last equality holds because $G_\# \Pi(\mu,\nu)$ is concentrated in Euclidean geodesic. Notice that I have found $\bar \pi \in \OptPlans(\mu,\nu)$ such that
 \begin{equation*}
     \int \di_{eu}^2(x,y) \de \bar \pi (x,y) < \int \di_{eu}^2(x,y) \de \Pi(\mu,\nu),
 \end{equation*}
 this contradicts the definition of $\Pi$.
 \end{pr}

In order to deduce the main result of this section I only have to prove the approximation property stated in Proposition \ref{prop:stabilitylast}, and to this aim I need to preliminary state and prove the following proposition. Let me also point out that this result can be proven using general theorems (see for example Theorem 10.27 in \cite{villani2008} or Theorem 1.3.1 in \cite{Figalli2010OptimalTA}), anyway I prefer to present a proof that uses only cyclical monotonicity arguments, similar to the ones explained previously.

\begin{prop}\label{prop:strconvexnorm}
 Let $\mathsf N:\R^N \to \R^+$ be a smooth norm, such that $\mathsf N^2:\R^N  \to \R^+$ is $k$-convex for some $k>0$. Calling $\di: \R^N \times \R^N \to \R^+$ the associated distance and given $\mu,\nu \in \ProbTwo(\R^N)$ with $\mu\ll\Leb^N$, there exists a unique $\pi\in \OptPlans(\mu,\nu)$ and it is induced by a map.
\end{prop}

\begin{pr}
According to Remark \ref{rem:uniquenessthroughmap}, it is sufficient to prove that every $\pi\in \OptPlans(\mu,\nu)$ it is induced by a map. To this aim, fix $\pi\in \OptPlans(\mu,\nu)$ and call $\Gamma$ the $\pi$-full measure, $\di^2$-cyclically monotone set, provided by Proposition \ref{prop:cmonotonicity}. Assume by contradiction that $\pi$ is not induced by a map, denote by $\{\pi_x\}_{x\in X}$ the disintegration kernel with respect to the projection map $\p_1$, then $\pi_x$ is not a delta measure for a $\mu$-positive set of $x$. Therefore there exists a compact set $A\subset \R^N$ with $\mu(A)>0$, such that $\pi_x$ is not a delta measure for every $x\in A$. Consequently consider 
\begin{equation*}
    \eta := \int_A \pi_x \times \pi_x \de \mu,
\end{equation*}
which is a positive measure on $\R^N \times \R^N$. Moreover $\eta$ is not concentrated on $\{(x,x):x\in\R^N\}$, thus there exists $(\bar x, \bar y)\in \supp (\eta)$ with $\bar x \ne \bar y$ and in particular $\eta (B_\delta (\bar x)\times B_\delta (\bar y))>0$ for every $\delta>0$. Now call $v=\bar y - \bar x$ and notice that, since $\mathsf N^2$ is smooth there exists $\bar \delta >0$ such that for every $z \in A$ it holds that
\begin{equation*}
    \modu{\frac{\partial}{\partial v} \mathsf N^2 (x-z)-\frac{\partial}{\partial v} \mathsf N^2 (\bar x-z)} < \frac k3 \di^2 (\bar y, \bar x)
\end{equation*}
for every $x \in B_{2 \bar \delta}(\bar x)$, and 
\begin{equation*}
    \modu{\frac{\partial}{\partial v} \mathsf N^2 (y-z)-\frac{\partial}{\partial v} \mathsf N^2 (\bar y-z)} < \frac k3 \di^2 (\bar y, \bar x)
\end{equation*}
for every $y \in B_{2\bar \delta}(\bar y)$. Moreover, since $\mathsf N^2$ is $k$-convex, for every $z \in A$ it holds that 
\begin{equation*}
    \frac{\partial}{\partial v} \mathsf N^2 (\bar y-z) \geq \frac{\partial}{\partial v} \mathsf N^2 (\bar x-z) + k \di^2 (\bar y, \bar x),
\end{equation*}
and consequently 
\begin{equation}\label{eq:partialder}
    \frac{\partial}{\partial v}\mathsf N^2 ( y-z) > \frac{\partial}{\partial v}\mathsf N^2 ( x-z) 
\end{equation}
for every $x \in B_{2 \bar \delta}(\bar x)$ and every $y \in B_{2\bar \delta}(\bar y)$. 
On the other hand, since $\eta (B_{\bar \delta} (\bar x)\times B_{\bar \delta} (\bar y))>0$, the set 
\begin{equation*}
    A^{\bar \delta} = \{ z \in \R^N : \pi_z (B_{\bar \delta} (\bar x))>0 \text{ and }\pi_z (B_{\bar \delta} (\bar y))>0\}.
\end{equation*}
has positive $\mu$-measure and thus it has positive $\Leb^N$-measure. Let $\bar z$ be the density point of $A^{\bar \delta}$, then in a neighborhood $\bar z$ there exists $z$ such that $z, z+ \epsilon v \in A^{\bar \delta}$ for some $0<\epsilon < \frac{\bar \delta}{\norm{v}}$. Consequently, it is possible to find $x\in B_{\bar \delta} (\bar x)$ and $y \in B_{\bar \delta} (\bar y)$, such that 
\begin{equation*}
    (z+ \epsilon v, x), (z,y) \in \Gamma.
\end{equation*}
Then it holds that 
\begin{align*}
    \di ^2(z,x) + \di^2 (z + \epsilon v,y) &= \mathsf N^2(x-z) +\mathsf N^2 (y-(z+\epsilon v))\\
    & =\mathsf N^2 (x-(z+\epsilon v)) + \int_0^\epsilon \frac{\partial}{\partial v}\mathsf N^2 (x-sv-z)) \de s\\
    & \quad +\mathsf N^2 (y-z) - \int_0^\epsilon \frac{\partial}{\partial v} \mathsf N^2 (y-sv-z) \de s \\
    & < \di ^2(z+ \epsilon v,x) + \di^2 (z ,y),
\end{align*}
where the last passage follows from \eqref{eq:partialder}. This last inequality contradicts the $\di^2$-cyclical monotonicity of $\Gamma$, concluding the proof.
\end{pr}

Having a consistent geodesic flow and an associated plan selection, it only remains to apply Proposition \ref{prop:stabilitylast} and deduce the main result. In order to do so, I introduce a sequence $(\di_n)_{n\in \setN}$ of distances on $\R^N$ by requiring the following three properties:
\begin{itemize}
\item for every n, $\di_n$ is induced by a smooth norm $\mathsf N$, such that $\mathsf N^2$ is $k$-convex for some $k>0$ and satisfies condition \eqref{eq:weirdcondition},
    \item $\di_n^2$ converges to $\di^2$ uniformly on compact sets,
    \item $n (\di^2_n - \di^2)$ converges to $\di^2_{eu}$ uniformly on compact sets, and $n (\di^2_n - \di^2)\leq2\di^2_{eu}$ for every $n$.
\end{itemize}
It is easy to see that such a sequence exists. Now, fixed a pair of absolutely continuous measures $\mu,\nu\in \Prob_{ac}(\R^N)$, Proposition \ref{prop:strconvexnorm} ensures that for every $n$ there exists a unique transport plan $\pi_n$ between $\mu$ and $\nu$, with respect to the cost $c(x,y)= \di_n^2(x,y)$. Let me now prove that it is possible to apply Proposition \ref{prop:stabilitylast}.

\begin{prop}
 The maps $G$ and $\Pi$ and the sequences $(\di_n)$ and $(\pi_n)$ I introduced satisfy the assumptions of Proposition \ref{prop:stabilitylast} with $K=0$.
\end{prop}

\begin{pr}
Condition 2 is easily satisfied, in fact since $\di_n$ is induced by a strictly convex norm the only geodesics in $(\R^N,\pi_n)$ are the Euclidean ones. Then, because $\pi_n$ is unique and Proposition \ref{prop:CDinR^n} holds, it is clear that the entropy functional is convex along $G_\#\pi_n$, with respect to the distance $\di_n$. 
Let me now prove condition 1. Notice that $\pi_n \in \Gamma(\mu,\nu)$ for every $n$, therefore the sequence $(\pi_n)$ is tight and Prokhorov theorem ensures the existence of $\pi\in \Gamma(\mu,\nu)$ such that, up to the extraction of a subsequence, $\pi_n\weakto \pi$. I am now going to prove that $\pi \in \Pi_2(\mu,\nu)$. Observe that $\pi_n$ is an optimal transport plan for the distance $\di_n$ and thus 
\begin{equation*}
    \int \di_n^2 \de \pi_n \leq \int \di_n^2 \de \tilde \pi \qquad \forall \tilde \pi \in \Gamma(\mu,\nu),
\end{equation*}
therefore for every compact set $C\subset \R^N$ it holds
\begin{equation*}
    \int_C \di_n^2 \de \pi_n \leq \int \di_n^2 \de \tilde \pi \qquad \forall \tilde \pi \in \Gamma(\mu,\nu).
\end{equation*}
It is then possible to pass to the limit as $n\to\infty$, using the uniform convergence for the left hand side and the dominated convergence (ensured by \eqref{eq:weirdcondition}) at the right hand side, obtaining 
\begin{equation*}
    \int_C \di^2 \de \pi \leq \int \di^2 \de \tilde \pi \qquad \forall \tilde \pi \in \Gamma(\mu,\nu).
\end{equation*}
Since this last equation holds for every compact set $C \subset \R^N$, it is possible to conclude that 
\begin{equation*}
    \int \di^2 \de \pi \leq \int \di^2 \de \tilde \pi \qquad \forall \tilde \pi \in \Gamma(\mu,\nu),
\end{equation*}
in particular $\pi \in \Pi_1(\mu,\nu)$. Using once more the minimizing property of $\pi_n$, follows that
\begin{equation*}
    \int \di^2 \de \tilde \pi + \int (\di_n^2- \di^2) \de \pi_n \leq \int \di_n^2 \de \pi_n \leq \int \di_n^2 \de \tilde \pi = \int \di^2 \de \tilde \pi + \int (\di_n^2- \di^2) \de \tilde \pi \qquad \forall \tilde \pi \in \Pi_1(\mu,\nu),
\end{equation*}
consequently it holds that
\begin{equation*}
    \int n(\di_n^2- \di^2) \de \pi_n \leq \int n(\di_n^2- \di^2) \de \tilde \pi \qquad \forall \tilde \pi \in \Pi_1(\mu,\nu),
\end{equation*}
and proceeding as before I can infer that
\begin{equation*}
    \int \di^2_{eu} \de \pi \leq \int \di^2_{eu} \de \tilde \pi \qquad \forall \tilde \pi \in \Pi_1(\mu,\nu).
\end{equation*}
In particular $\pi\in \Pi_2(\mu,\nu)$ and this concludes the proof, considering the definition of the map $\Pi$.
\end{pr} 

\noindent Finally, the combination of this last result with Proposition \ref{prop:stabilitylast} allows me to conclude the final result of this article.

\begin{corollary}
The metric measure space $(\R^N,\di,\Leb^N)$ is a very strict $\CD(0,\infty)$ space and consequently it is weakly essentially non-branching.
\end{corollary} 

\noindent {\scshape Aknowlegments} : This article contains part of the work I did for my master thesis, that was supervised by Luigi Ambrosio and Karl-Theodor Sturm.

\bibliography{verystrictCD}

\begin{thebibliography}{10}

\bibitem{AmbrosioNotes2}
L.~Ambrosio.
\newblock Lecture notes on optimal transport problem.
\newblock Euro Summer School "Mathematical aspects of evolving interfaces",
  2000.

\bibitem{ambrosio2013user}
L.~Ambrosio and N.~Gigli.
\newblock A user’s guide to optimal transport.
\newblock In {\em Modelling and optimisation of flows on networks}, pages
  1--155. Springer, 2013.

\bibitem{Ambrosio_2015}
L.~Ambrosio, N.~Gigli, A.~Mondino, and T.~Rajala.
\newblock Riemannian {R}icci curvature lower bounds in metric measure spaces
  with $\sigma $-finite measure.
\newblock {\em Transactions of the American Mathematical Society},
  367(7):4661–4701, 2015.

\bibitem{Ambrosio_2013}
L.~Ambrosio, N.~Gigli, and G.~Savaré.
\newblock Calculus and heat flow in metric measure spaces and applications to
  spaces with {R}icci bounds from below.
\newblock {\em Inventiones mathematicae}, 195(2):289–391, 2013.

\bibitem{Ambrosio_2014}
L.~Ambrosio, N.~Gigli, and G.~Savaré.
\newblock Metric measure spaces with {R}iemannian {R}icci curvature bounded
  from below.
\newblock {\em Duke Mathematical Journal}, 163(7):1405–1490, 2014.

\bibitem{ambkirpra03}
L.~Ambrosio, B.~Kirchheim, and A.~Pratelli.
\newblock Existence of optimal transport maps for crystalline norms.
\newblock {\em Duke Mathematical Journal}, 125(2):207--241, 2004.

\bibitem{Figalli2010OptimalTA}
A.~Figalli.
\newblock {\em Optimal Transportation and Action-Minimizing Measures}.
\newblock Edizioni della Normale. SNS, 2010.

\bibitem{gigli11}
N.~Gigli.
\newblock Optimal maps in non branching spaces with {R}icci curvature bounded
  from below.
\newblock {\em Geometric and Functional Analysis}, 22:990--999, 2011.

\bibitem{Gigli_2015}
N.~Gigli, A.~Mondino, and G.~Savaré.
\newblock Convergence of pointed non-compact metric measure spaces and
  stability of {R}icci curvature bounds and heat flows.
\newblock {\em Proceedings of the London Mathematical Society}, 111:1071--1129,
  2015.

\bibitem{lottvillani}
J.~Lott and C.~Villani.
\newblock Ricci curvature for metric-measure spaces via optimal transport.
\newblock {\em Annals of Mathematics}, 169:903--991, 2009.

\bibitem{MM-Example}
M.~Magnabosco.
\newblock Example of an highly branching {CD} space.
\newblock {\em arXiv preprint}, 2021.

\bibitem{rajala2013failure}
T.~Rajala.
\newblock Failure of the local-to-global property for {$\CD (K, N)$} spaces.
\newblock {\em Ann. Sc. Norm. Super. Pisa Cl. Sci.}, 15:45--68, 2016.

\bibitem{rajalasturm}
T.~Rajala and K.-T. Sturm.
\newblock Non-branching geodesics and optimal maps in strong {$\CD(K,\infty)$}
  spaces.
\newblock {\em Calculus of Variations and Partial Differential Equations},
  50:831--846, 2014.

\bibitem{santambrogio2015optimal}
F.~Santambrogio.
\newblock {\em Optimal Transport for Applied Mathematicians: Calculus of
  Variations, PDEs, and Modeling}.
\newblock Progress in Nonlinear Differential Equations and Their Applications.
  Springer International Publishing, 2015.

\bibitem{schultz2017existence}
T.~Schultz.
\newblock Existence of optimal transport maps in very strict {$\CD(K,\infty)$}
  spaces.
\newblock {\em Calculus of Variations and Partial Differential Equations}, 57,
  2018.

\bibitem{Schultz2019EquivalentDO}
T.~Schultz.
\newblock Equivalent definitions of very strict {$\CD(K,N)$} spaces.
\newblock {\em arXiv preprint}, 2019.

\bibitem{Schultz2019OnOO}
T.~Schultz.
\newblock On one-dimensionality of metric measure spaces.
\newblock {\em Proc. Amer. Math. Soc.}, 149:383--396, 2020.

\bibitem{sturm2006}
K.-T. Sturm.
\newblock On the geometry of metric measure spaces.
\newblock {\em Acta Math.}, 196(1):65--131, 2006.

\bibitem{sturm2006ii}
K.-T. Sturm.
\newblock On the geometry of metric measure spaces. {II}.
\newblock {\em Acta Math.}, 196(1):133--177, 2006.

\bibitem{villani2003}
C.~Villani.
\newblock {\em Topics in Optimal Transportation}.
\newblock Graduate studies in mathematics. American Mathematical Society, 2003.

\bibitem{villani2008}
C.~Villani.
\newblock {\em Optimal transport -- Old and new}.
\newblock Grundlehren der mathematischen Wissenschaften. Springer, 2008.

\end{thebibliography}

\bibliographystyle{abbrv} 

\end{document}